\documentclass[a4book, 10pt]{amsart}
\title[]{Flat Surfaces with singularities\\
 in Euclidean 3-space}

\usepackage[dvips]{graphicx} 
\usepackage{verbatim,enumerate}
\usepackage{amssymb}
\usepackage{amsthm}
\theoremstyle{plain}
 \newtheorem{theorem}{Theorem}[section]
 \newtheorem*{theorem*}{Theorem}

 \newtheorem*{lemma*}{Lemma}
 \newtheorem{proposition}[theorem]{Proposition}
 \newtheorem{fact}[theorem]{Fact}
 \newtheorem*{fact*}{Fact}

 \newtheorem{lemma}[theorem]{Lemma}
 \newtheorem{corollary}[theorem]{Corollary}
\theoremstyle{remark}
 \newtheorem{definition}[theorem]{Definition}
 \newtheorem{remark}[theorem]{Remark}
 \newtheorem*{remark*}{Remark}
 \newtheorem*{acknowledgements}{Acknowledgements}
 \newtheorem{example}[theorem]{Example}
\numberwithin{equation}{section}

\newcommand{\Z}{\boldsymbol{Z}}
\newcommand{\R}{\boldsymbol{R}}
\renewcommand{\phi}{\varphi}
\renewcommand{\epsilon}{\varepsilon}

\newcommand{\op}{\operatorname}

\newcommand{\mc}[1]{{\mathcal #1}}
\newcommand{\mb}[1]{{\mathbf #1}}
\newcommand{\pmt}[1]{{\begin{pmatrix} #1  \end{pmatrix}}}
\newcommand{\dy}{\displaystyle}

\author{Satoko~Murata}
\address[Murata]{%
   Department of Mathematics, Kyoto Tachibana High School,
    Momoyamacho Iga 50, Kyoto 612-8026, Japan
}
\email{murata@tachibana-hs.jp}

\author{Masaaki Umehara}
\address[Umehara]{%
   Department of Mathematics, Graduate School of Science,
   Osaka University,
   Toyonaka, Osaka 560-0043,
   Japan
}
\email{umehara@math.sci.osaka-u.ac.jp}

\thanks{
The second author was partly supported by the Grant-in-Aid for 
Scientific Research (B), Japan Society for the Promotion of Science.}

\begin{document}
\begin{abstract}
It is classically known that complete flat 
(that is, zero Gaussian curvature)
surfaces in Euclidean 3-space
$\R^3$ are cylinders over space curves.
This implies that the study of global behaviour 
of flat surfaces requires the study of
singular points as well.
If a flat surface $f$ admits singularities but
its Gauss map $\nu$ is globally defined on the surface
and can be smoothly extended
across the singular set, $f$ is called a {\it frontal}.
In addition, if the pair $(f,\nu)$ defines an immersion 
into $\R^3\times S^2$, $f$ is called a {\it front}. 
A front $f$ is called {\it flat} if the Gauss map
degenerates everywhere.
The parallel surfaces and the
{\it caustic} (i.e.~focal surface) 
of a flat front $f$ are also flat fronts.
In this paper, we generalize the classical notion of
completeness to flat fronts,
and give a representation formula for a
flat front which has  a non-empty compact singular set and
whose ends are all immersed and complete.
As an application, we show that such a flat 
front has properly embedded ends
if and only if its Gauss map image is a convex curve.
Moreover, we show the existence of at least
four singular points other than cuspidal edges on such a flat front with
embedded ends, which is a variant of the classical four vertex theorem for
convex plane curves. 
\end{abstract}
\maketitle

%%%%%%%%%%%%%%%%%%%%%%%%%%%%%%%%%%%%%%%%%%%%%%%%%%%%%
\section*{Introduction}
In this paper, we shall investigate the global 
behaviour of flat surfaces
with singularities in Euclidean 3-space $\R^3$.
In fact, for the study of global properties of flat surfaces,
considering only immersions is too restrictive, as
is already clear from the
classical fact (Fact 0.1) below.
 
Let $M^2$ be a smooth $2$-manifold and 
$f:M^2\to \R^3$ a $C^\infty$-map.
A point $p\in M^2$ is called {\it regular} if
$f$ is an immersion on a sufficiently small neighborhood
of $p$, and is called {\it singular} if it is not regular. 
If $f$ is an immersion and has zero Gaussian curvature,
it is called a (regular) flat surface.
It is classically known that regular
flat surfaces have open dense subsets 
which consist of a cylinder, cone, 
or tangential developable of 
space curves. 
Moreover, the following result is also well-known:

\begin{fact}\label{fact:I}
A complete regular flat surface is a plane or a cylinder.
\end{fact}

\noindent
The first proof of this fact came as a corollary of Hartman
and Nirenberg \cite{HN}. After that 
Stoker \cite{S} and Massey \cite{M} gave elementary proofs.

To extend flat surfaces to a larger class that allows singularities,
we shall now  give the following definitions:
A $C^\infty$-map $f:M^2\to \R^3$ is called a {\it frontal} if
there exists a smooth 
unit vector field $\nu$ along $f$ such that
$\nu$ is perpendicular to $df(TM^2)$.
By parallel displacement, $\nu$ can be
considered as a map into the unit sphere $S^2$.
In this case, $\nu$ is called the {\it Gauss map}
of the frontal $f$.
Moreover, if the map
$$
L:=(f,\nu):M^2\to \R^3\times S^2
$$
gives an immersion, $f$ is called a {\it front} or a {\it wave front}.
Using the canonical inner product,
we identify the unit tangent bundle
$
\R^3\times S^2=T_1\R^3
$
with the unit cotangent bundle $T_1^*\R^3$, which has
the canonical contact structure.
When $f$ is a front, $L$ gives a Legendrian immersion.
Hence, wave fronts are considered as projections of  
Legendrian immersions.

For a frontal  $f:M^2\to \R^3$ and a real number $\delta$,
a new frontal  $f_\delta:M^2\to \R^3$ is given by
\begin{equation}\label{eq:parallel}
f_\delta:=f+\delta \nu 
\end{equation}
which is called a {\it parallel surface} of $f$.
A frontal $f$ is called {\it flat} if the Gauss map $\nu:M^2\to S^2$
degenerates.
Since the members of the parallel family $\{f_\delta\}$ have the same Gauss map
$\nu$ in common, {\it parallel surfaces of a flat front are also flat}. 
Moreover, {\it the caustics  $($i.e. focal surfaces$)$ of  flat fronts
are also flat} (Section 1).

Flat surfaces might not have globally defined 
normal vector fields, in general.
A $C^\infty$-map $f:M^2\to \R^3$ is called a {\it 
p-frontal} (resp.~{\it p-front}) if for each $p\in M^2$, there exists a
neighborhood $V$ such that the restriction $f|_V$
of $f$ to $V$ gives a  frontal  (resp.  front).
That is, p-front (i.e. projective front) 
is locally a front, but the unit normal
vector field $\nu$ might not extend globally.
A p-frontal is called {\it co-orientable} if
it is a front (that is, $\nu$ is globally defined on $M^2$).
The class of flat p-fronts is a canonical generalization of that of 
regular flat surfaces.
In fact, the existence of a real analytic  flat M\"obius strip is known
(\cite{W}),
which is an example of a flat p-front without singular points.

We will now give two definitions for 
the completeness of flat fronts. 

\begin{definition}\label{def:complete}
A flat p-front $f:M^2\to \R^3$ 
with first fundamental form $ds^2=df\cdot df$
is called {\it complete} (cf. \cite[p150]{KRUY}) if
there exists a symmetric covariant tensor $T$ on $M^2$
with compact support
such that $ds^2+T$ gives a complete metric
on $M^2$. 
If $f$ is complete and the singular set is non-empty, 
then the singular set must be compact, by definition.
Moreover, this new definition of completeness for flat fronts 
is a generalization of the classical one, namely, 
if $f$ is complete and has no singular points, then
completeness coincides with 
the classical notion of
completeness in Riemannian Geometry. 

On the other hand, $f$ is called {\it weakly complete}
if the sum of the first fundamental form and the
third fundamental form (called the 
{\it lift metric} of the p-front $f$)
\begin{equation}\label{eq:lift}
ds^2_\#:=df\cdot df+d\nu\cdot d\nu
\end{equation}
gives a complete Riemannian metric on $M^2$, which
is the pull-back of the canonical metric on 
$T_1\R^3$ by $L=(f,\nu)$.
\end{definition}

Obviously, {\it completeness implies weak completeness.}
Completeness is convenient to describe
the asymptotic behaviour of ends, but
is not preserved when lifting to
the universal cover.
On the other hand, the universal cover
of a weakly complete flat front is also weakly complete. 
Anyhow, using these definitions, we will prove our deepest
result (Theorem A below):

\medskip
\noindent
{\bf Theorem A.}
{\it Let $M^2$ be a connected 2-manifold, and
$f:M^2\to \R^3$ a complete flat p-front
whose singular set is non-empty $($namely, the
singular set of $f$ is compact and all ends of $f$ are
immersed and complete$)$. Then it has no umbilic points and is a front.
Moreover, $M^2$ is diffeomorphic to a circular cylinder.}

\medskip
In particular, if a complete flat p-front admits an
umbilic point, it must be a cylinder.
We remark that Theorem~A has two surprising consequences:
\begin{itemize}
\item A complete flat p-front is a front.
In particular,  it never contains a M\"obius strip as a subset.
\item Singular points and umbilics never appear on the
same complete flat front.
\end{itemize}

Moreover, as a consequence of Theorem~A, we get  
the following representation formula for complete flat fronts:

\medskip
\noindent
{\bf Theorem B.}
{\it Let $\hat \xi:S^1\to S^2$ be a regular curve 
without inflection points,
and $\alpha=a(t) dt$ a 1-form on $S^1=\R/2\pi \Z$ satisfying the
period condition
$$
\left(\int_{S^1}\hat\xi \alpha=\right)
\int_{0}^{2\pi}\hat\xi(t) a(t) dt=0,
$$
where $\hat\xi \alpha$ is considered as a
$2\pi$-periodic  $\R^3$-valued 1-form 
on $\R$.
Then 
\begin{equation}\label{eq:B}
f(t,v)=\hat\sigma(t)+v\hat \xi(t),\quad
\hat\sigma(t):=\int_{0}^t \hat\xi(t) a(t)dt
\end{equation}
gives a complete flat front whose image of the singular set
coincides with that of $\hat\sigma$. 
Conversely, any complete flat fronts with non-empty
singular set are given in this manner. 
}

\medskip
A somewhat similar representation theorem for 
(regular) flat tori in the 3-sphere has been given
by Kitagawa \cite{Ki},
that is, a pair of two closed regular spherical curves 
satisfying a certain compatibility condition
correspond to a  flat torus in $S^3$.
As a consequence of Theorem~B,
 outside of a ball of sufficiently large radius
centered at the origin,
the image of a complete
flat front splits
into two components, called the
{\it ends} of the front.
We shall prove the following:

\medskip
\noindent
{\bf Theorem C.}
{\it 
An end of a complete flat front $f$ with a singular point 
is properly embedded if and only if the other end is as well.
Moreover, the ends of $f$ are properly embedded
if and only if the image of the Gauss map of $f$ is a convex curve.}

\medskip
It was shown in \cite{KUY} that a given complete flat front in 
hyperbolic 3-space $H^3$ has embedded ends
if and only if the sum of the topological 
degrees of the two hyperbolic Gauss maps 
of the surface coincides with the number of ends.
(In contrast to our case, there are complete flat fronts
in $H^3$ with arbitrary genus and with many ends.)
The above theorem is an analogue for $\R^3$,
though our method is quite different from that of
\cite{KUY}.
As an application of Theorem~B, Theorem~C
and the technique of disconjugate operators (cf. \cite{gmo} and \cite{tu2}),
we prove in this paper the following

\medskip
\noindent
{\bf Theorem D.}
{\it Let $f:M^2\to \R^3$ be a complete flat front
with non-empty singular set.
Suppose that all ends of $f$ are embedded.
Then there exist at least four 
singular points of $f$ on $M^2$
which are not cuspidal edges.
}

\medskip
The proof of Theorem~D is an analogue of 
that of the classical four vertex theorem for convex
plane curves. 
The embeddedness of ends is truly required in the above statement.
In fact, there exists a complete flat front in $\R^3$
which admits only cuspidal edges. (See Example \ref{ex:twist}.)
We show in Section 4 (see Corollary~\ref{cor:sing}) that
complete flat fronts which admit only
cuspidal edges and swallowtails are generic. 
For such a generic front with embedded ends,
Theorem~D implies the existence of at least four swallowtails.

In Section 1, we give fundamental properties of flat fronts.
In Section 2, we define admissibility of developable
frontal and give
a representation formula
for admissible developable frontals.
In Section 3, we extend the representation formula
to p-fronts. As an application, we show the existence of
umbilic points on a developable M\"obius strip.
In Section 4, Theorem~A and Theorem~B are proved.
In Section 5, we shall prove Theorem~C and Theorem~D. 

Finally, we mention several related works:
\cite{GM} for flat fronts in $S^3$;
\cite{SUY}, \cite{FSUY} for behaviour of
the Gaussian curvature near cuspidal edges, 
swallowtails, and cuspidal cross caps;
\cite{UY}, \cite{KY}, \cite{F}, \cite{FRUYY},\cite{LY} for 
singularity of spacelike maximal surfaces 
in $\R^3_1$ and CMC-1 surfaces in de Sitter 3-space; 
\cite{M} for improper affine spheres with singularities; 
 \cite{IM}
for generic singularities of fronts of constant Gaussian curvature 
in $\R^3$.

\begin{acknowledgements}
The authors thank Wayne Rossman, Kotaro Yamada and the referees 
for a careful reading of the first draft 
and for valuable comments. 
Their comments were especially helpful for making the revision.
The authors also thank Osamu Saeki 
for valuable comments.
\end{acknowledgements}

\section{Flat surfaces as wave fronts}\label{sec:prelim}
We fix a 2-manifold $M^2$ throughout this section.
First we recall the following fundamental property of
fronts: 

\begin{lemma}\label{lem:1-1}
Let $M^2$ be a 2-manifold and $p\in M^2$
a singular point of a front $f:M^2\to \R^3$.
Then for any sufficiently 
small non-zero real number $\delta$,
the point $p$ is a regular point of 
the parallel surfaces $f_\delta$. 
\end{lemma}

\begin{proof}
Since $p$ is a singular point, the differential of the map
$(df)_p:T_pM^2\to T_{f(p)}\R^3$ is of rank less than $2$.
When $(df)_p$ is of rank zero, then 
the Gauss map $\nu$ of $f$  is an immersion 
at $p$, and the parallel surface $f_\delta:=f+\delta\nu$ is
an immersion at $p$ for all $\delta \ne 0$.
So it is sufficient to consider the case that 
$(df)_p$ has a 1-dimensional kernel.
Then, we can take a local coordinate $(U;u,v)$
such that $f_u$ vanishes at $p$, where
$f_u,f_v$ are partial derivatives of the 
$\R^3$-valued function $f$
with respect to the parameters $u,v$.
Since $f_u(p)=0$, we have  
\begin{equation}\label{eq:1-1}
0=f_u(p) \cdot \nu_v(p) =-f_{uv}(p)\cdot \nu(p) 
=f_v (p)\cdot \nu_u(p),
\end{equation}
where \lq$\cdot$\rq\ denotes the canonical inner product on
Euclidean 3-space $\R^3$.
Since $f$ is a front and
$(df)_p$ has 1-dimensional kernel,
$\nu_u(p), f_v(p)$ are both non-vanishing.
Thus \eqref{eq:1-1} yields that $\nu_u,f_v$ are
linearly independent at $p$, 
and the assertion follows immediately.
\end{proof}

Let $f:M^2\to \R^3$ be an immersion
with globally defined unit normal vector field $\nu$,
and $\lambda_1,\lambda_2$ two (distinct) 
principal curvature functions.
At umbilic points the two functions coincide,
and they are continuous functions on $M^2$ and
are $C^\infty$ at non-umbilic points.
The following fact is well known:

\begin{fact}\label{fact:1-2}
The principal curvature functions 
$\lambda_{1,\delta},\lambda_{2,\delta}$ of 
the regular set of each parallel surface
$f_\delta:=f+\delta\nu$ \rm ($\delta\in \R$) satisfy
\begin{equation}\label{eq:change}
\lambda_j=\lambda_{j,\delta}(1-\delta \lambda_j)
\qquad (j=1,2).
\end{equation}
Moreover, the singular set $S(f_\delta)$ of $f_\delta$ 
is given by 
\begin{equation}\label{eq:Sdelta}
S(f_\delta):=\bigcup_{j=1,2}
\{p\in M^2\,;\, \delta\lambda_j(p)=1\}.
\end{equation}
\end{fact}

Now let
$
f:M^2\to \R^3
$
be a front. We fix a point $p\in M^2$ arbitrarily.
By Lemma \ref{lem:1-1}, there exist a neighborhood
$U$ of $p$ and a real number $c$ such that
the parallel surface $f_c$ gives an immersion on $U$.
A point $p$ is called an {\it umbilic point} of $f$ if
it is the umbilic point of one such $f_c$.
By definition, the set of umbilic points $\mc U_f$ 
is common in its parallel family.

\begin{lemma}\label{lem:wein}
Let
$
f:M^2\to \R^3
$
be a front and $\nu$ the unit normal vector field.
For each non-umbilic point $p$, there exists a 
local coordinate $(U;u,v)$ of $M^2$ centered at $p$
such that
$\nu_u$ and $\nu_v$
are proportional to $f_u$ and $f_v$ on $U$
respectively. 
In particular, the pair of 
vectors $\{f_u,\nu_u\}$ $($resp. $\{f_v,\nu_v\})$
do not vanish at the same time.
\end{lemma}

\begin{remark}\label{rmk:coordinate}
Such a coordinate $(U;u,v)$ is called
a {\it curvature line coordinate} of 
the front $f$. 
If $f$ is an immersion,
the principal curvature functions 
$\lambda_1,\lambda_2$
satisfy 
\begin{equation}\label{weingarten}
\nu_u=-\lambda_1 f_u,\qquad
\nu_v=-\lambda_2 f_v,
\end{equation}
which is called the Weingarten formula.
The above lemma is a generalization
of it for fronts.
\end{remark}

\begin{proof}
By Lemma \ref{lem:1-1}, there exists 
a parallel surface $f_\delta$, which is regular at $p$.
Since $p$ is not an umbilic point, 
there exists a curvature line coordinate 
$(U;u,v)$ satisfying 
$\nu_u=-\lambda_{1,\delta} (f_\delta)_u,\,\,
\nu_v=-\lambda_{2,\delta} (f_\delta)_v$.
Since $f=f_\delta-\delta\nu$,
we get the assertion.
\end{proof}

\begin{definition}
The direction $v$ in $T_pM^2$ such that
$df(v),d\nu(v)$ are linearly dependent is called 
a {\it principal direction} of a front $f$.
\end{definition}

It can be easily checked that 
the two  principal directions are common in the parallel family of a front.
 
Let $f:M^2\to \R^3$
be a front with the unit normal vector field $\nu$.
We fix $p\in M^2\setminus \mc U_f$ arbitrarily.
By Lemma \ref{lem:wein}, 
we can take a curvature line coordinate $(U;u,v)$ 
containing $p$.
Then $\{f_u,\nu_u\}$ and $\{f_v,\nu_v\}$ are
linearly dependent respectively.
So we define maps 
$
\Lambda_j:M^2\setminus \mc U_f\to P^1(\R)\,\,(j=1,2)
$
by
$$
\Lambda_1=[-\nu_u:f_u],\qquad \Lambda_2=[-\nu_v:f_v],
$$
where $[-\nu_u:f_u]$ and $[-\nu_v:f_v]$
mean the proportional ratio of $\{-\nu_u,f_u\}$
and $\{-\nu_v,f_v\}$ respectively as elements of the real 
projective line $P^1(\R)$.
When $f$ is an immersion, then
we can write $\nu_u=-\lambda_1f_u$, $\nu_v=-\lambda_2f_v$, 
and we have $\Lambda_j=[-\lambda_j,1]\,\,(j=1,2)$.
If $f_u=\nu_u=0$ at $p$, then $f$ cannot be a front at $p$.
So $\Lambda_1$ and $\Lambda_2$ are both 
well-defined $C^\infty$-maps,
which are called the {\it principal curvature maps}.
They are canonical generalizations of the principal
curvature functions.

\begin{proposition}
The principal curvature maps
$\Lambda_1,\Lambda_2$ can be continuously
extended to the entirety of $M^2$.
\end{proposition}

\begin{proof}
We fix $p\in M^2$ arbitrarily.
By Lemma \ref{lem:1-1}, 
there exists a parallel surface $f_\delta$, which is regular 
on a neighborhood $U$ of $p$.
Then  there are two principal curvature functions
$\lambda_{1,\delta},\lambda_{2,\delta}$.
We set
$$
\hat \Lambda_j=
[-\lambda_{j,\delta}:1+\delta \lambda_{j,\delta}]
\qquad (j=1,2).
$$
Since $\lambda_{1,\delta},\lambda_{2,\delta}$ can be
taken to be continuous on $U$, the maps
$\hat\Lambda_1,\hat\Lambda_2$
are both continuous on $U$.
Since $p$ is arbitrary, $\hat\Lambda_1,\hat\Lambda_2$ are
both continuous functions on $M^2$.
On the other hand, by Fact \ref{fact:1-2},
$\hat\Lambda_1,\hat\Lambda_2$ coincide with
the principal curvature maps of $f$
defined on $M^2\setminus \mc U_f$, so the assertion is
proved.
\end{proof}

The following assertion is now obvious from 
the definition of $\Lambda_1$ and $\Lambda_2$.

\begin{proposition}\label{prop:1-3}
Let $f:M^2\to \R^3$ be a front, and
$\Lambda_1,\Lambda_2$ the principal curvature maps.
Then  a point $p\in M^2$ is 
an umbilic point if
$\Lambda_1(p)=\Lambda_2(p)$.
On the other hand,
$p$
is a 
singular point if and only if
$\Lambda_1(p)=[1:0]$
or $\Lambda_2(p)=[1:0]$.
\end{proposition}

As defined in the introduction, a front $f:M^2\to \R^3$ is 
called {\it flat} if its Gauss map degenerates everywhere.
Since the Gauss map is common in the same parallel family,
we get the following assertions.

\begin{proposition}\label{prop:gaussian0}
Let $f:M^2\to \R^3$ be a front
and $\nu$ the unit normal vector field.
Then $f$ is flat if and only if, 
for each point $p$, there exist a 
real number $\delta$ and a neighborhood 
$U$ of $p$
such that
the parallel surface $f_\delta$
is an immersion on $U$ and
its Gaussian curvature vanishes
identically on $U$.
\end{proposition}

\begin{corollary}\label{lem:1-5}
Let $f:M^2\to \R^3$ be a front.
Suppose that the regular set $R(f)$ of $f$ is dense in $M^2$.
Then $f$ is flat if and only if $f$ has zero Gaussian curvature
on $R(f)$.
\end{corollary}

\begin{proposition}\label{prop:umbilics}
Let $p$ be an umbilic point of a flat front
$f:M^2\to \R^3$, then $p$ is a regular point of $f$.
\end{proposition}

\begin{proof}
By the definition of umbilic points for fronts, 
both of $\nu_u$ and $\nu_v$ vanish at $p$.
Since $(f,\nu)$ is an immersion, $f$ must be
an immersion. Thus $p$ is a regular point of $f$.
\end{proof}

\begin{corollary}\label{lem:1-4}
Let $f:M^2\to \R^3$ be a flat front and $p\in M^2$
a singular point. Then $p$ is a regular point of $f_\delta$
for all $\delta\ne 0$.
\end{corollary}

\begin{proof}
Since $p$ is a singular point, it is not an umbilic point.
Thus the assertion follows from Lemma \ref{lem:1-1} and 
(1.3), since we may set one of the principal curvature
functions to be vanishing identically.
\end{proof}

Now we fix a flat front $f:M^2\to \R^3$ with a
unit normal vector field $\nu$.
Then we may assume that
$
\Lambda_2=[0:1]
$
identically.
(That is, $v$-lines correspond to the asymptotic direction.)
So we may set
$$
\Lambda:=\Lambda_1:M^2\to P^1(\R),
$$
which is called  the {\it principal curvature map} 
of the flat front $f$.
Moreover, the function
$$
\rho:=\frac{\psi_2}{\psi_1}:M^2\to \R\cup\{\infty\}
\qquad (\Lambda=[-\psi_1:\psi_2])
$$
is called the {\it  curvature radius function}.
If $f$ is an immersion, $\rho$ coincides with the usual
curvature radius function.
The restriction of $\rho$ to $M^2\setminus \mc U_f$
is a real valued $C^\infty$-function. 

\begin{definition}\label{def:asymptotic}
Let $f:M^2\to \R^3$ be a flat front, and
$J$ an open interval. 
A curve $\gamma:J\to M^2$ is called an {\it asymptotic curve}
if $\hat\gamma''(t)\cdot \nu$ vanishes identically, where
$\hat\gamma=f \circ \gamma(t)$ and $\prime=d/dt$.
\end{definition}

Now, we fix a non-umbilic point $p\in M^2$
and take a  special {\it curvature line coordinate}
$(U; u,v)$ centered at $p$, such that
\begin{enumerate}
\item $u$-curves are non-zero principal curvature lines,
and $v$-curves are asymptotic curves.
\item $|f_{v}|=1$, that is $v$ is the arc-length parameter of
$v$-curves. 
\end{enumerate}
The existence of $(U; u,v)$ for flat immersion is 
well-known (See \cite{M} or \cite[5-8]{C}).
Since $\nu_v=0$, $|(f_\delta)_v|=|f_v|=1$ holds for parallel 
surface, the existence of the coordinate
for flat front is also shown.
The following assertion can be proved directly.

\begin{proposition}\label{prop:1-9}
On the curvature line coordinate $(U;u,v)$
as above, 
the derivative of the Gauss map 
$\nu_u$ never vanishes
and satisfies
$\nu_v=0,\,\, f_{vv}=0$.
\end{proposition}

\begin{proof}
Since
$$
f_v\cdot \nu_u=(f_v \cdot \nu)_u-f_{uv} \cdot \nu
=-f_{uv} \cdot \nu=
(f_u \cdot \nu)_v-f_{u} \cdot \nu_v=0,
$$
$\{f_v,\nu,\nu_u\}$
is a local orthogonal frame field on $U$.
Then $|f_v|=1$ implies that $f_{vv}$ is perpendicular
to $f_v$. Obviously, $f_v$ is orthogonal to $\nu$.
Thus, to prove $f_{vv}=0$, it is sufficient to
show that $f_{vv}$ is perpendicular to $\nu_u$.
In fact, we have
$$
f_{vv}\cdot \nu_u=(f_{v}\cdot \nu_u)_v-f_{v}\cdot \nu_{uv}
=-f_{v}\cdot (\nu_{v})_u=0,
$$
which proves the assertion.
\end{proof}

In the proof of Fact 0.1 in \cite{M},
the following important lemma is given (cf. \cite[5-8]{C}):
\begin{fact}\label{lem:Massey} $($Massey's lemma$)$
Let $f:M^2\to \R^3$ be a flat immersion.
Let $J$ be an open interval and
$\gamma(v)$ $(v\in J)$
be an asymptotic curve on $M^2$ passing through
a non-umbilic point $p$.
Suppose that $v$ is the arc-length parameter on it.
Then $\hat\gamma:=f\circ \gamma$ is contained in a straight line
and the restriction $\rho(v)$ of the curvature radius function satisfies
the equality $\rho''(v)=0$.
Furthermore, the closure of $\gamma(J)$ 
does not meet the umbilic set $\mc U_f$. 
\end{fact}
The following generalization of Massey's lemma holds for
flat fronts:
\begin{lemma} \label{lem:gmassey} $($Generalized Massey's lemma$)$
The above Massey's lemma $($Fact \ref{lem:Massey}$)$ 
holds for any flat front $f:M^2\setminus \mathcal U_f\to \R^3$.
In particular, the curvature radius function $\rho$ is
a smooth function on 
$M^2\setminus \mathcal U_f$ 
which satisfies $\rho''(v)=0$ 
along each asymptotic curve $\gamma(v)$
with $|\hat\gamma'(v)|=1$.
\end{lemma}

\begin{proof}
Since the asymptotic curve is a zero curvature line,
it is also an asymptotic curve of the parallel surface $f_\delta$.
For a given  $p\in M^2$,
there exists a real number $\delta$ 
and a neighborhood $U$ of $p$ such that
$f_\delta$ is an immersion on $U$.
Then the image of the asymptotic curve $f_\delta\circ \gamma
=\hat\gamma(v)+\delta\nu(\gamma(v))$ is contained in 
a straight line. Since the Gauss map
$\nu$ is constant along $\gamma$, the image of 
$\hat \gamma=f\circ \gamma$ is also contained in  a straight line. 

On the other hand, by Fact \ref{lem:Massey},
the principal curvature radius
function $\rho_\delta$ (of $f_\delta$)
satisfies
$\rho_\delta''(v)=0$
whenever $\gamma(v)\in U$
by  Lemma \ref{lem:Massey}.
Then, we have
$$
0= \rho''_\delta=(\rho-\delta)''=\rho'',
$$
which proves the first part.
(The property that $\prime=d/dv$
is the arclength parameter of $f_\delta$
does not depends on $\delta$.)
Next we suppose that $\gamma(J)$
accumulates at $q\in \mc U_f$.
Then there exists a sequence $\{p_n\}$
in $\gamma(J)$ converging to $q$.
Since $q$ is a regular point of $f$,
there exists a neighborhood $U$ of $q$
such that $f$ is an immersion on $U$.
Without loss of generality, we may 
assume that $\{p_n\}$ lies in $U$.
Since $f|_U$ is a flat immersion, by Fact \ref{lem:Massey}, 
$\{p_n\}$ cannot accumulate at 
the umbilic point $q$, which is a contradiction.
\end{proof}

Let $\rho$ be the curvature 
radius function of $f$.
Then the map
$$
C_f:=f+\rho \nu:M^2\setminus \mc U_f\to \R^3
$$
is called the {\it caustic} or the {\it focal surface} 
of the front $f$.
The following assertion can be proved by a straight-forward calculation
using curvature line coordinates as in Proposition \ref{prop:1-9}:

\begin{theorem}\label{thm:caustics}
The caustic  
$C_f:M^2\setminus \mc U_f\to \R^3$ 
of a flat front $f$
is also a flat front.
Moreover, the caustic has no umbilics.
\end{theorem}

Later, we shall reprove the above assertion as an application
of our representation formula for flat fronts in Section 2.
(See Proposition~\ref{prop:caustics}.)
The caustics of flat fronts in the hyperbolic 3-space $H^3$ 
or the 
3-sphere $S^3$ are also flat.
Roitman \cite{R} gave a representation formula for caustics of
flat fronts in $H^3$. (See also \cite{KRSUY}.)

%%%%%%%%%%%%%%%%%%%%%%%%%%%%%%%%%%%%%%%%%%%%%%%%%%%%%%
\section{The structure of developable frontals } 

In this section, we shall investigate a weaker class
of flat surfaces:
Let $M^2$ be a 2-manifold.
Recall that a $C^\infty$-map $f:M^2\to \R^3$
is called a {\it frontal} if there exists a
unit vector field $\nu:M^2\to S^2$ such that
$\nu(p)$ is perpendicular to $df(T_pM^2)$ for
each $p\in M^2$.
($\nu$ is called the {\it unit normal vector field}
of $f$.)
In addition, if the pair 
$(f,\nu):M^2\to \R^3\times S^2$ 
defines an immersion, $f$ is called a {\it front}. 

A frontal  $f:M^2\to \R^3$ is called 
a {\it ruled frontal}, if for each $p\in M^2$,
there exists a neighborhood $U$ of $p$
and a vector field $\xi$ on $U$ such that
\begin{enumerate}
\item $df_q(\xi_q)$ is a unit vector for all $q\in U$,
\item the image of each integral curve of $\xi$ by $f$
is a part of straight line in $\R^3$. 
\end{enumerate}
This $\xi$ is called a {\it $($local$)$ unit asymptotic vector field} 
of $f$.
Moreover, if the Gauss map $\nu:M^2\to S^2$ of $f$ 
is constant on each integral curve of $\xi$,
$f$ is called a {\it developable frontal}.
In the introduction, {\it flat frontal} is defined.
We note that {\it developable frontals are examples of 
flat frontals.}
In this case, the integral curves associated with
$\xi$ are called the {\it asymptotic lines}.
(When $f$ is a front, asymptotic lines coincide with
the asymptotic curve as in Definition~\ref{def:asymptotic}.)

A vector field $X$ on a manifold $M$ is called
{\it complete} if each integral curve $\gamma(t)$ of $X$ is 
defined for $t\in \R$. 
If $X$ is complete, it generates a 1-parameter
group of transformation $\phi_t:M\to M$ ($t\in \R$).

On the other hand,
we denote by $P(TM)$ the projective tangent bundle,
whose fibre is the projective space $P(T_pM)$ of the tangent 
vector space $T_pM$. A section of $P(TM)$ is called a
{\it projective vector field}. Let 
 $\xi$ be a projective vector field on $M$ which is not a 
globally defined vector field.
Then there exists double covering 
$\hat M$ of $M$ such that $\xi$ can be lifted up as a
$C^\infty$-vector field $\hat \xi$ on $\hat M$.
In this case, the integral curves of $\hat \xi$ correspond to
those of $\xi$, but the orientation of the curves as a flow
does not have global meaning on $M$ in general.

\begin{definition}\label{def:asm-complete}
A ruled frontal $f:M^2\to \R^3$ is called
 {\it asymptotically complete} or
{\it $a$-complete} if the unit asymptotic vector field $\xi$
is determined as a complete projective vector 
field, that is, all of the asymptotic lines of $f$ 
are complete, namely, each image of them by $f$ consists 
of the entirety of a straight line.
\end{definition}

\begin{proposition}\label{ftodev}
Let $f:M^2\to \R^3$ be a flat front.
If $f$ has no umbilic points, then
$f$ is developable.
\end{proposition}

\begin{proof}
We take a local curvature line coordinate $(U;u,v)$
of $f$ as in Proposition~\ref{prop:1-9} such that
$v$-curves correspond to the zero curvature lines
satisfying $|f_v|=1$.
Then 
$
\xi_U=\partial/\partial v
$
gives the asymptotic direction, and
$\nu_v=0$ holds.  
\end{proof}

\begin{definition}\label{def:sigma}
A ruled frontal $f:M^2\to \R^3$ is called
{\it admissible} if 
it is $a$-complete and
there exists an embedding
$\sigma:J\to M^2$ such that
\begin{enumerate}
\item $J$ is an open interval of $\R$ or $J=S^1$,
\item $\sigma'$ 
is transversal to the unit asymptotic vector field
$\xi$,
\item $\sigma$ meets each asymptotic line exactly once.
\end{enumerate}
This curve $\sigma$ is called a {\it generator} 
of an admissible developable frontal $f$,
which is not uniquely determined in general
(see Remark~\ref{rem:canonical}).
\end{definition}

\begin{example}\label{ex:dev}
Let 
$
c(t):J\to \R^3$ ($J=\R$ or $S^1$)
be a smooth map and 
$\xi(t)\in T_{c(t)}\R^3$ a
smooth vector field along $c$
such that
\begin{equation}\label{eq:flat}
\det(c'(t),\xi(t),\xi'(t))=0 \qquad (t\in J).
\end{equation}
Then
$$
F_{c,\xi}(t,v):=c(t)+v\xi(t)
\qquad (t\in J)
$$
gives a typical example of a developable frontal, if
there is a unit vector field $\nu(t)$
of $F_{c,\xi}$ along $c$ such that $c'(t)\cdot \nu(t)=
\xi(t)\cdot \nu (t)=0$
for all $t\in J$. For example, when
$c'(t)$ and $\xi(t)$ are linearly independent,
$F_{c,\xi}$ is a developable frontal,
since $\nu=c'\times \xi/|c'\times \xi|$ is
a unit normal vector field of $F_{c,\xi}$.
In this situation, $F_{c,\xi}$ is
admissible, since 
the $t$-axis in $(t,v)$-plane meets 
all of the $v$-curves exactly once transversally.
Later, we shall show that an admissible
flat frontal is congruent to one of
$F_{c,\xi}$ (see Proposition~\ref{prop:dev_str}).
Let $c(t)$ be a regular space curve. 
Then $
\xi(t):=c'(t)
$
satisfies \eqref{eq:flat}, and
the induced map
$$
F_{c,\xi}(t,v):=c(t)+vc'(t)
\qquad (t\in J)
$$
is called the {\it tangential developable},
whose singular set is $c(J)$.
\end{example}

Now we return to the general situation.
Let $f:M^2\to \R^3$ be an $a$-complete ruled frontal.
If the unit asymptotic vector field
$\xi$ can be globally defined as a single-valued vector field on $M^2$,
$f$ is called {\it $a$-orientable} or {\it asymptotically orientable}.
There are non-$a$-orientable developable fronts. 
(See Proposition \ref{prop:non-ori} in the next section.)
In this section, we shall mainly consider
$a$-orientable developable frontals.

\begin{proposition}\label{prop:dev_str}
Let $f:M^2\to \R^3$ be an $a$-complete
and  $a$-orientable ruled frontal,
and $\xi$ the unit asymptotic vector field on $M^2$.
Suppose that $f$ is admissible and
$
\sigma:J\to M^2
$
is a generator $($cf. Definition~\ref{def:sigma}$)$, 
where $J=\R$ or $S^1$.
Then there exists a diffeomorphism
$
\Phi:J\times \R\to M^2
$
such that
$$
f\circ \Phi(t,v)=
f\circ \sigma(t)+v
df(\xi_{\sigma(t)})\qquad (|df(\xi_{\sigma(t)})|=1).
$$
In particular, $M^2$ is orientable.
\end{proposition}

\begin{proof}
Since $f$ is $a$-complete, each integral curve
$\phi_v(q)$ of $\xi$ passing through $q$
is defined for all $v\in \R$.
In particular, it induces a 
1-parameter group of transformations
$$
\phi_v:M^2\to M^2\qquad (v\in \R). 
$$
Now we define a smooth map $\Phi$ by
$$
\Phi:J\times \R 
\ni (t,v)\mapsto \phi_v(\sigma(t))
\in M^2.
$$
Since
$\sigma'$ is transversal to $\xi$,
so is $d\phi_v(\sigma'(t))$, which
implies that $\Phi$ is locally a diffeomorphism.
The bijectivity of $\Phi$ follows from
the fact that $\sigma(J)$ meets all 
asymptotic lines exactly once. 
Thus $\Phi$ is a diffeomorphism.
By definition, for each $t\in J$
$$
f\circ \Phi(t,v)=f(\phi_v(\sigma(t)))\qquad (v\in \R)
$$
gives the straight line
passing through $f(\sigma(t))$. Moreover, since
$$
\left|\frac{\partial f\circ 
\Phi(t,v)}{\partial v}\right|=|df(\xi)|=1,
$$
the parametrization $v\mapsto f\circ \Phi(t,v)$
coincides with the arc-length parameter of the
straight line passing through $f(\sigma(t))$
with the initial velocity vector $df(\xi_{\sigma(t)})$.
Thus we get
$$
f\circ \Phi(t,v)=
f\circ \sigma(t)+v
df(\xi_{\sigma(t)}),
$$
which proves the assertion.
\end{proof}

\begin{remark}\label{rem:canonical}
The above coordinate $(t,v)$ of $M^2$
is called a {\it canonical coordinate of $M^2$
with respect to a generator $\sigma$}.
In this coordinate, each $t$-curve
$J\ni t\mapsto \Phi(t,v)\in M^2$ 
(for  fixed $v\in \R$)
gives also a generator.
So, for an arbitrarily  given point $p\in M^2$,
we can take a canonical coordinate $(t,v)$
such that $\Phi(0,0)=p$.
\end{remark}

Now we fix an admissible
flat developable with a canonical coordinate
$$
f(t,v)=f\circ\sigma(t)+v
df(\xi_{\sigma(t)})\qquad
(t\in J,\,\, v\in \R),
$$
such that $|df(\xi_{p})|=1$ for all $p\in M^2$.
Then we get a smooth curve 
$$
\hat\xi:J\ni t
\mapsto df(\xi_{\sigma(t)})\in S^2,
$$
called the {\it asymptotic Gaussian curve},
which gives a parametrization of
the image of the asymptotic Gauss map
$
M^2\ni p\mapsto 
df(\xi_{p})
\in S^2.
$

On the other hand, a smooth curve
$$
\hat\nu:J\ni t
\mapsto \nu_{\sigma(t)}\in S^2
$$
is called the {\it Gaussian curve} of $f$,
which gives a parametrization of
the image of the  Gauss map
$\nu:M^2\to S^2$. \\

The common parameter $t$ of two curves $\hat \xi(t)$ 
and $\hat \nu(t)$ depends on the choice
of a generator of $f$, but their
spherical images $\hat \xi(J)$ and $\hat \nu(J)$
do not depend on the choice
of generator.

\begin{proposition}\label{prop:duality}
Let $f:M^2\to \R^3$ be an admissible $a$-orientable 
developable frontal.
Then the Gaussian curve $\hat \nu(t)$ 
and the asymptotic Gaussian curve $\hat \xi(t)$
satisfy the following relations
$$
\hat\nu\cdot \hat \xi=\hat\nu'\cdot \hat\xi=\hat\nu\cdot \hat \xi'=0,
$$
that is, $\hat\nu$ and $\hat\xi$ are dual frontals in $S^2$,
where $'=d/dt$. 
$($See the appendix.$)$
Moreover, when $f$ is a front, $\hat \nu'(t)$ 
vanishes if and only if $\nu^{-1}(\hat \nu(t))$ is an umbilic set.
In particular, if $f$ has no umbilics, then $\hat \nu(t)$ is
a regular curve, and $\hat\xi(t)$ is a front $($as 
a spherical curve$)$. 
\end{proposition}

\begin{proof}
The first assertion
$\hat\xi\cdot \hat\nu=0$ obviously holds.
We fix $p\in M^2$ arbitrarily.
As mentioned in Remark \ref{rem:canonical},
we can take a canonical coordinate 
$$
\Phi:J\times \R\ni (t,v)\mapsto \Phi(t,v)\in M^2
$$
such that $\Phi(0,0)=p$ and 
$|d(f\circ \Phi)(\partial/\partial v)|=1$.
Then
\begin{align*}
-\hat \xi' \cdot \hat \nu&=
\hat \xi \cdot \hat \nu'=
f_{v}(t,0)\cdot \nu_t(t,0)\\
&=-f_{vt}(t,0)\cdot \nu(t,0)
=f_{t}(t,0)
\cdot \nu_v(t,0)=0.
\end{align*}
Next, we suppose that $f$ is a front and
$\hat \nu'(t)=0$, then it implies 
$\nu_t(0,0)=0$, and so $d\nu(0,0)=0$,
that is, $p$ is an umbilic point of $f$.
\end{proof}

Let $f:M^2\to \R^3$ be an admissible $a$-orientable 
developable frontal.
In the same notations as in the proof of
Proposition~\ref{prop:duality},
we have the expression 
$$
f\circ \Phi(t,v)=\hat\sigma(t)+v \hat\xi(t)
\qquad (\,(t,v)\in J\times \R\,),
$$
where
$
\hat\sigma(t)=f\circ \Phi(t,0).
$
Since $\hat\sigma'(t)$ is perpendicular to $\hat\nu$,
we can write
$$
\hat\sigma'=a \hat \xi +b \hat \eta,
\qquad \hat \eta:=\hat \xi \times \hat \nu,
$$
where \lq$\times$\rq\ is the vector product in $\R^3$,
and $a,b\in C^\infty(J)$. 
We define 1-forms $\alpha,\beta$ on $J$ by
$$
\alpha:= a(t) dt,\quad 
\beta:= b(t) dt \qquad (t\in J). 
$$
Since
$
d\hat \sigma=\alpha \hat \xi+\beta \hat \eta
$
holds, the quadruple $(\alpha,\beta,\hat \xi,\hat \nu)$
is independent of the choice of the parameter $t$ on 
$J$ as a 1-dimensional manifold. 
It induces a $C^\infty$-map
$$
f(t,v)=v \hat\xi+
\int_0^t \left(\alpha \hat \xi+\beta \hat \eta\right),
\qquad \hat\eta:=\hat\xi\times \hat\nu,
$$
which coincides with $f\circ \Phi$ up to parallel translation.
We prove the following:

\begin{theorem}\label{thm:w-rep}
$($A representation formula for 
developable frontals$)$ \\
Let $\alpha,\beta$ be 1-forms on an open 
interval $J (\subset \R)$.
Let $\hat \nu,\hat \xi:J\to S^2$ be 
mutually dual frontals.
Then a $C^\infty$-map defined by
\begin{align}\label{eq:rep_open}
f(t,v)&:=
v \hat\xi(t)+\hat\sigma(t),\\
\nonumber
\hat\sigma(t)&:=
\int_0^t \left(\alpha\hat \xi +\beta \hat \eta\right),
\quad \hat\eta:=\hat\xi\times \hat\nu
\end{align}
gives an $a$-orientable
admissible developable frontal on $J\times \R$.
The singular set of $f$ is
the zeros  of $\mu_1(t) v+b(t)$,
where 
$$
\hat \xi'(t)=\mu_1(t) \hat \eta(t),\qquad
\beta=b(t)dt.
$$
Moreover, $f$ is a front 
if and only if $\hat \nu'(t)\ne  0$ or
\begin{equation}\label{front-cond}
\hat\xi'(t)=\hat\nu'(t)=0 \,\,\mbox{and}\,\, b(t)\ne 0
\end{equation}
holds for each $t\in J$.
$($In particular, if $\hat \nu$ is a regular curve, 
then $f$ is a front.$)$
Conversely, any admissible $a$-orientable developable frontals
defined on $J\times \R$ are given in this manner.
\end{theorem}

\begin{remark}\label{rem:front}
As stated in the theorem, $\hat \nu'(t)\ne 0 
\,\,(t\in J)$ is
a sufficient condition that $f$ is a front.
On the other hand, if $f$ is a front,
then the condition \eqref{front-cond} at $t=t_0$ implies that $f$ 
is an immersion at $\gamma(t_0)$ and the asymptotic line
passing through $\gamma(t_0)$ consists of umbilic points.
For example, if we set
$$
\alpha=0,\quad \beta=dt,\quad
\hat \xi(t)=(1,0,0),\quad \hat \nu(t)=(0,0,1),
$$
then the data $(\alpha,\beta,\hat \xi, \hat\nu)$
satisfies \eqref{front-cond},
and the corresponding surface is the plane 
$f(t,v)=(t,v,0)$.
In this case, $\hat \nu'(t)$ vanishes identically
although $f$ is an immersion.
\end{remark}

\begin{remark}\label{rem:variants}
Quadruples 
$$
(\alpha,\beta,\hat \xi, \hat\nu),\quad 
(\alpha,\beta,-\hat \xi, -\hat\nu),
\quad 
(-\alpha,-\beta,\hat \xi, \hat\nu)
$$
give the same developable frontals up to congruence.
Moreover, in our construction, the quadruple 
$$
(\alpha+d\phi,\,\, \beta+\phi \mu_1 dt,
\,\, \hat \xi,\,\, \hat\nu)
$$
corresponds to a flat developable $f+\phi\hat\xi$ 
which is congruent to $f$.
\end{remark}

\begin{proof}
We have already seen that
any $a$-complete flat frontals
defined on $J\times \R$
are given in this manner.
So it is sufficient to show that $f(t,v)$ given 
by \eqref{eq:rep_open} is an (orientable) $a$-complete developable frontal.
It can be easily checked that $\hat\nu$ 
gives a unit normal vector field of $f$.
Since it does not depend on $v$, 
the flatness follows immediately.
Since $v\in \R$ is arbitrary, $f$ is obviously $a$-complete
and admissible. 
(The curve $t\mapsto (t,0)$ is a generator.)

Finally, we shall examine the condition that
$f$ is a front.
Suppose that $f$ is a frontal and fix a point $p=(t_0,v_0)\in M^2$
arbitrarily. Then 
$\hat \nu'(t_0)\ne 0$ holds if and only if
$f$ is a front at $p$, and $p$
is not an umbilic point.
So it is sufficient to show that \eqref{front-cond}
holds if and only if $f$ is an immersion at $p$
and $\nu'(t_0)=0$.
In fact, if $f$ is a front and $p$ is an umbilical point,
we have $\hat \nu'(t_0)=0$, that is, $d\nu$ vanishes at $(t_0,v)$
for all $v\in \R$.
Thus $f(t,v)$ must be an immersion at $(t_0,v)$ ($v\in \R$).
On the other hand, we have that 
\begin{align}\label{eq:immersion}
f_t\wedge f_v&=(a(t)\hat \xi(t)+b(t)\hat \eta(t)+v \hat \xi'(t))\wedge
\hat\xi(t)\\
&=(b(t)+v \mu_1(t))\hat\eta(t)\wedge \hat\xi(t),\nonumber
\end{align}
where $\hat \xi'(t)=\mu_1(t)\hat \eta(t)$. 
If $\mu_1(t_0)\ne 0$, $f(t_0,v)$ is not an immersion at $v=-b(t)/\mu_1(t)$,
a contradiction. Thus $\mu_1(t_0)= 0$, and then we have $b(t_0)\ne 0$
since $f_t\wedge f_v\ne 0$.
Conversely, \eqref{front-cond} and \eqref{eq:immersion}
implies that $f$ is an immersion at $(t_0,v)$ for all $v\in \R$.
\end{proof}

Next, we consider the case 
$
J=S^1=\R/2\pi\Z.
$
The following corollary is immediate from
the above theorem.

\begin{corollary}\label{cor:rep-1closed}
$($A representation formula for 
developable frontals with closed generator$)$
Let $\alpha,\beta$  be 1-forms on $S^1$,
and $\hat \nu,\hat \xi:S^1\to S^2$ mutually
dual frontals.
Suppose that
\begin{equation}\label{eq:rep-1closed}
\int_{S^1} \left(\alpha\hat \xi+\beta\hat \eta\right)=0,
\quad \hat\eta:=\hat\xi\times \hat\nu.
\end{equation}
Then a $C^\infty$-map defined by
\begin{equation}\label{eq:rep-2}
f(t,v):=\hat\sigma(t)+
v \hat\xi(t),
\qquad \hat\sigma(t):=
\int_0^t \left(\alpha\hat \xi+\beta\hat \eta\right)
\end{equation}
is an admissible developable frontal on $S^1\times \R$.
Conversely, any  admissible orientable developable frontals
defined on $S^1\times \R$ are given in this manner.
\end{corollary}

The simplest developable front with closed generator 
is the circular cylinder 
$$
f(u,v)=(\cos t,\sin t,v),
$$
which corresponds to the data
 $(0,dt,\mathbf e_3,f(t,0))$, where $\mathbf 
e_3=(0,0,1)$.
We give here another example:

\begin{example}\label{ex:rectifying}
(Rectifying developables)
Let $J=\R$ or $S^1$, and
let $c:J\to \R^3$ be a regular space curve
satisfying $c''(t)\ne 0$ for $t\in J$.
We may assume that $t$ is the arc-length parameter.
Then the envelope $f$ of the rectifying plane of $c(t)$
is called a {\it rectifying developable}
associated with $c$. It is well-known that
$f$ has an expression
$
f(t,v)=c(t)+v\xi(t)
$,
where
$$
\xi(t):=\frac{1}{\sqrt{\kappa(t)^2+\tau(t)^2}}
(\kappa(t) \mb e(t)+\tau(t) \mb b(t)),
$$
is called the {\it Darboux vector field}.
(Here $\mb e(t)=c'(t)$,
$\mb b(t)$ is the unit bi-normal vector,
$\kappa(t), \tau(t)$ are the curvature function and the
torsion function of $c(t)$ respectively.)
When $\kappa(t)\ne 0$, it corresponds
to the data  
$\left(\frac{\kappa(t)dt}{\sqrt{\kappa^2+\tau^2}},
\frac{-\tau(t)dt}{\sqrt{\kappa^2+\tau^2}}
,\xi ,\mb n\right),
$ where $\mb n=\mb n(t)$ is the 
unit principal normal vector of $c(t)$.
\end{example}

We return now to the general situation.

\begin{proposition}\label{prop:p-period}
Let $f(t,v)$ be an admissible developable front
as in Corollary~\ref{cor:rep-1closed}.
Suppose that $f$ has no umbilics.
Then the non-zero 
principal curvature lines of $f$ are
characterized by a graph $v=v(t)$ in $tv$-plane
satisfying $v'(t)+a(t)=0$.
In particular, the non-zero 
principal curvature lines are
all periodic if
$
\int_{S^1}a(t)dt=0.
$
\end{proposition}

\begin{proof}
In the $tv$-plane,
each non-zero principal curvature line
can be expressed by
$w(t):=(t,v(t))$,
since it is transversal to the asymptotic
lines (i.e.~$v$-lines).
Since the principal directions are common in
the parallel family and our assertion is local,
we may assume that $w(t)$ lies in the regular
set of $f$.
Since $f_t-af_v$ is proportional to $\hat \eta$,
it is perpendicular to
the asymptotic direction $\hat \xi$.
Thus, $w(t)$ gives the non-zero principal direction 
if and only if
$$
\frac{d(f\circ w(t))}{dt}=f_t+v'f_v,
$$ 
is proportional to $f_t-af_v$, which is equivalent to
the condition $v'+a=0$.
Now the assertion follows immediately.
\end{proof}

\begin{remark}\label{rem:infinite}
The linear map
$$
C^\infty(S^1,\R^2)\ni
(a,b)\mapsto \int_{S^1}
\left(
 a(t)\hat \xi(t)+b(t)\hat \eta(t),
a(t)\right)dt
\in \R^4
$$
has infinite dimensional kernel
corresponding to
the set of $a$-complete developable fronts whose curvature lines are 
periodic. 
Thus there are infinitely many 
examples whose curvature lines are 
periodic. For example,
a cone over a locally convex spherical curve
given in Example \ref{cone} in Section 4 
is a complete flat front whose 
non-zero curvature lines are all closed, since $a(t)$
vanishes identically.
\end{remark}

Let $f:S^1\times J\to \R^3$ be an admissible 
developable frontal corresponding to
a quadruple $(a(t)dt,b(t)dt,\hat\xi(t),\hat\nu(t))$.
We now define two $C^\infty$-functions
$\mu_1(t)$ and $\mu_2(t)$ on $J$ by
\begin{equation}\label{eq:mus}
\hat\xi'=\mu_1 \hat\eta,\qquad \hat\nu'=\mu_2 \hat \eta,
\qquad \hat\eta:=\hat\xi\times \hat\nu
\qquad (\prime=d/dt).
\end{equation}
Then the singular set is equal to the 
zeros of $\mu_1(t) v+b(t)$
by Theorem~\ref{thm:w-rep}. 

\begin{definition}\label{def:linear}
A singular point 
$p=(t,v)$ is {\it linear} if and only if 
$\hat\nu'(t)=0$ (i.e. $\mu_1(t)=0$).
\end{definition}

If $p$ is a linear singular point,
the asymptotic line passing through $p$
consists of linear singular points.
The singular set of the cylinder over a cardioid
(cf. Remark \ref{rmk:cardioid})
is a typical example of linear singular points.

A singular point is called a 
{\em cuspidal edge\/} or {\em swallowtail\/} 
if it is locally diffeomorphic to 
\begin{equation}\label{eq:cuspidal-swallow}
 f_C(u,v):=(u^2,u^3,v), \qquad
 f_S(u,v):=(3u^4+u^2v,4u^3+2uv,v)
\end{equation}
at $(u,v)=(0,0)$, respectively.
These two types of singular points are 
the generic singularities of fronts. 
We recall here criteria 
given in \cite{KRSUY} as follows:
Let $f:(U;u,v)\to \R^3$ be 
a front (not necessarily flat).
We denote by $\nu(u,v)$, the unit normal vector
field of $f(u,v)$. Then there exists a smooth function 
$\lambda(u,v)$ defined on $U$ such that 
$f_u\times f_v=\lambda \nu$.
A point $p\in U$ is called {\it non-degenerate} 
 if and only if $d\lambda\ne 0$.
In this case, the implicit function theorem yields that
the singular set consists of a regular curve around $p$,
and such a  curve is called a {\it singular curve}.
We denote by $\sigma(t)$
the singular curve passing through $p$.
Since $\sigma$ is a regular curve on $U$,
the direction $\sigma'(t):=d \sigma(t)/dt$ is 
called the {\it singular direction}.
On the other hand, since $f(u,v)$ is not an immersion 
along $\sigma$,
there exists a non-vanishing vector field $Z(t)$
(along $\sigma$)
such that $df(Z(t))=0$.
The following criteria are known $($see \cite{KRSUY}.$)$
when $f$ is a front.
\begin{itemize}
\item[(a)] {\it $p=\sigma(t_0)$ 
is a cuspidal edge if and only if 
$p$ is a non-degenerate singular point
and $Z(t_0)$ is not proportional to $\sigma'(t_0)$,
that is, the determinant function $\phi(t):=\det(Z(t),\sigma'(t))$
does not vanish at $t=t_0$.}
\item[(b)] {\it $t=t_0$ is a swallowtail if and only if 
$p$ is a non-degenerate singular point such that
$\phi(t_0)=0$ and 
$\phi'(t_0)\ne 0$.}
\end{itemize}

\begin{proposition}\label{prop:singular}
Let $p=(t_0,v_0)$ be a singular point of an 
admissible developable front $f(t,v)$
as in Theorem~\ref{thm:w-rep}.
\begin{enumerate}
\item[{\rm (1)}] $f$ is non-degenerate at $p$ if $p$
is a non-linear singular point or a linear singular point
satisfying $b'(t_0)+v_0 \mu'_1(t_0)\ne 0$.
\item[{\rm (2)}] Suppose that $p$ is non-linear, namely
$\mu_1(t_0)\ne 0$, then
\begin{enumerate}
\item[{\rm (i)}] a singular point $p$ is locally diffeomorphic
to a cuspidal edge if and only if
$\mu_2(t_0)\ne 0$ and $a\ne (b/\mu_1)'$ at $t=t_0$,
\item[{\rm (ii)}] $p$ is locally diffeomorphic
to a swallowtail if and only if
$\mu_2(t_0)\ne 0$, $a(t)=(b(t)/\mu_1(t))'$ and
$a'(t)\ne (b(t)/\mu_1(t))''$ at $t=t_0$.
\end{enumerate}
\item[{\rm (3)}] Suppose that $p$ is linear, namely
$\mu_1(t_0)=0$, then $b(t_0)$ must vanish and
$p$ is locally diffeomorphic
to a cuspidal edge if and only if
$\mu_2(t_0)\ne 0$ and $v_0\mu'_1(t_0)+b'(t_0)\ne 0$.
Moreover, swallowtails never appear at $p$.
\end{enumerate}
\end{proposition}

\begin{remark}
In \cite{SUY}, the singular curvature $\kappa_s$
is defined on cuspidal edges. 
Since the Gaussian curvature of $f$ is
non-negative, $\kappa_s$
is non-positive as shown in \cite[Theorem~3.1]{SUY}.
Since $|\kappa_s|$ is equal to the 
tangential component of the second derivative of the
singular curve, the equality of 
$\kappa_s(p)\le 0$ holds if and only if
$p$ is a linear singular
point.
\end{remark}

\begin{proof}
Our conclusions follow  from the above criteria (a) and (b).
In fact, $f_u\times f_v=(v\mu_1+b)\nu$, and
$p=(t_0,v_0)$ is a singular point if and only if
$
v_0\mu_1(t_0)+b(t_0)=0.
$
So the singular point $p=(t_0,v_0)$ is 
linear i.e. $\mu_1(t_0)=0$, and then we have $b(t_0)=0$.
On the other hand, $p$ is non-degenerate if and only if
$\mu_1(t_0)\ne 0$ or $v\mu'_1+b'\ne 0$ at $t=t_0$.
In particular, non-linear singular points are always
non-degenerate, and we get (1).
Suppose $p$ is a non-degenerate singular point.
Then the singular curve $\sigma$ 
passing through $p=(t_0,v_0)$
is given by
$$
\sigma(t)=
\begin{cases}
(t,-b(t)/\mu_1(t)) & \mbox{(if $\mu_1(t_0)\ne 0$)},\\
(t_0,t)            & \mbox{(if $\mu_1(t_0)=0$)}.
\end{cases}
$$
Since $f_t-a f_v$ vanishes on this singular curve,
the null direction is $Z(t)=(1,-a(t))$.
So we have
$$
\det(\sigma',Z)=
\begin{cases}
\det\pmt{1 & 1 \\ -(b/\mu_1)'& -a}=-a+(b/\mu_1)' & (\mbox{if}\,\,
\mu_1(t_0)\ne 0),\\
\det\pmt{0 & 1 \\ 1& -a}=-1 & (\mbox{if}\,\,\mu_1(t_0)= 0),
\end{cases}
$$
where \lq$\det$' denotes the determinant of two column-vectors.
The assertions  (2) and (3) now follow immediately. 
When $\mu_1(t_0)=0$,
a non-degenerate singular point is always a cuspidal edge.
In particular, a swallowtail never appears.
\end{proof}

In the previous 
section (cf. Theorem~\ref{thm:caustics}), 
we pointed out that the caustics of flat 
fronts are again flat.
We give here a refinement of this.

\begin{proposition}\label{prop:caustics}
Let $f$ be a front associated with  
$(adt,bdt,\hat\xi,\hat\nu)$.
Then the quadruple 
$
(a dt,(b+\delta \mu_2)dt,\hat\xi,\hat\nu)
$ 
corresponds to the parallel front $f_\delta$
\,\,$(\delta\in \R)$,
where $\hat\nu'=\mu_2\hat\eta$.
Suppose that $f$ is a front and has no umbilics
$($namely $\mu_2(t)$ does not vanish 
and $\hat \nu$ is a regular spherical curve$)$.
Then the caustic $C_f$ corresponds to the
quadruple 
$$
\left(\frac{a \mu_2 +(b/\mu_2)'
\mu_1 }{\sqrt{(\mu_1)^2+(\mu_2)^2}}dt,
\frac{-a\mu_1+(b/\mu_2)'\mu_2}{\sqrt{(\mu_1)^2
+(\mu_2)^2}}dt
,c_{\xi},\hat\eta\right),
$$
where $\prime=d/dt$, $\hat\xi'=\mu_1\hat\eta$ and
$
c_{\xi}:=\frac{\mu_2\hat \xi-\mu_1 \hat\nu}
{\sqrt{(\mu_1)^2+(\mu_2)^2}}.
$
In particular, the caustic $C_f$ of an admissible 
developable front is  also admissible.
\end{proposition}

\begin{proof}
The first assertion is obvious.
The non-zero principal 
curvature radius function of $f(t,v)$ is given by
$$
\rho(t)=-\frac{b(t)+v\mu_1(t)}{\mu_2(t)}.
$$
Since 
$C_f(t)=f(t,v)+\rho(t) \hat\nu(t)$, 
we have
$$
(C_f)_v=\hat \xi-\frac{\mu_1}{\mu_2}\hat \nu, \qquad
(C_f)_t=a \hat \xi-\left(\frac{b+\mu_1\hat\nu}{\mu_2} \right)'\hat \nu. 
$$
Thus $\hat \eta$ is perpendicular to both $(C_f)_v$ and $(C_f)_t$,
and gives a unit normal vector of $C_f$.
Since $\hat \eta$ depends only on the parameter $t$,
the Gauss map of $C_f$ degenerate everywhere. 
Moreover, since $\hat \eta'=-\mu_1\hat \xi-\mu_2\hat \nu$
and $\mu_2\ne 0$, $C_f$ gives an $a$-complete flat front.
On the other hand, since
$$
C_f=\left(\hat \sigma-\frac{b}{\mu_2}\hat \nu\right) 
+v\mu_2\sqrt{(\mu_1)^2+(\mu_2)^2}\,c_\xi,
$$
$c_\xi(t)$ gives the asymptotic Gaussian curve of $C_f$.
Furthermore, by a straight-forward calculation, we have
$$
\left(\hat \sigma-\frac{b}{\mu_2}\hat \nu\right)'
=\frac{a \mu_2 +(b/\mu_2)'\mu_1 }
{\sqrt{(\mu_1)^2+(\mu_2)^2}}c_\xi
+\frac{-a\mu_1+(b/\mu_2)'\mu_2}{\sqrt{(\mu_1)^2
+(\mu_2)^2}}
\mb e_f,
$$
where
$$
\mb e_f(t):=c_\xi(t) \times \hat \nu(t)
=-\frac{\mu_1(t) \hat \xi(t)+\mu_2(t) \hat \nu(t)}
{\sqrt{\mu_1(t)^2+\mu_2(t)^2}}.
$$
This proves the assertion.
\end{proof}

\begin{remark}\label{rem:duality}
By Proposition~A.1. in the appendix, $c_{\xi}(t)$ as above
coincides with the caustic of the asymptotic Gaussian curve 
$\hat \xi(t)$ of $f$.
\end{remark}

%%%%%%%%%%%%%%%%%%%%%%%%%%%%%%%%%%%%%%%%%%%%%%%%%%%%%%%%%%%%%%
\section{Non-orientable or non-co-orientable
developable surfaces.}

In this section, we shall generalize
the representation formula for admissible developables
to those which are not orientable or not admitting 
a globally defined unit normal vector field.

\begin{definition}
Let $M^2$ be a 2-manifold.
As defined in the introduction,
a $C^\infty$-map $f:M^2\to \R^3$ is called
a {\it p-frontal}  $($resp. {\it p-front}$)$
if for each $p\in M^2$,
there exists a neighborhood $U$ of $p$ such that
the restriction $f|_U$ gives a frontal 
$($resp. front$)$.
A p-frontal $f$ is called {\it co-orientable} if
it is a frontal, that is,
there exists a smooth unit normal vector field
globally defined on $M^2$. 
\end{definition}

If $f:M^2\to \R^3$ is non-co-orientable p-frontal,
there exists a double covering $\pi:\hat M^2\to M^2$
such that $f\circ \pi$ is a frontal.

\begin{definition}
A p-frontal is called {\it developable}
if $f$ is a developable  frontal or
 $f\circ \pi$ is a developable frontal. 
Similarly, {\it $a$-completeness of a p-frontal} is defined 
using the double covering.
A developable p-frontal is called {\it $a$-orientable}
if and only if there is a globally defined unit asymptotic 
vector field.
Like as the case of $frontal$ (cf. Definition~\ref{def:sigma}),
we define the admissibility for p-frontals:
A developable 
p-frontal $f:M^2\to \R^3$ is called {\it admissible} if
it is $a$-complete, and there is a curve 
(called a {\it generator})
embedded in $M^2$ satisfying the
three conditions in Definition~\ref{def:sigma}.
\end{definition}

For example, a real analytically immersed flat M\"obius strip
 (cf. \cite{W})
is non-orientable and non-co-orientable at the same time.
The following lemma is a key to proving 
the main results in this section.

\begin{lemma}\label{lem:transversal}
$($A criterion of admissibility$)$
Let $f:M^2\to \R^3$ be 
an  
$a$-complete
developable p-frontal.
Suppose that $M^2$ is connected and 
there exists an embedded closed curve
$$
\sigma:S^1=\R/2\pi \Z\to M^2,
$$
whose velocity vector $\sigma'(t)$ is
transversal to the asymptotic direction
of $f$ at $\sigma(t)$ for each $t\in S^1$,
and meets each asymptotic line in at most one point.
Then $f$ is admissible and $\sigma$
is a generator. 
\end{lemma}

\begin{proof}
Firstly, we suppose that
$f$ is $a$-orientable.
Then we can take a unit asymptotic vector field $\xi$
defined on $M^2$.
Since $f$ is $a$-complete,
each integral curve of vector field $\xi$ is 
defined on $\R$. 
In particular, $\xi$ generates a
global 1-parameter group of transformations
$$
\phi_v:M^2\to M^2\qquad (v\in \R)
$$
such that  
$v\mapsto f(\phi_v(p))$
is a straight line in $\R^3$.\\
Let $C$ be the image of the curve $\sigma$.
It is sufficient to show that $C$ meets 
the asymptotic line $v\mapsto \phi_v(p)$ for each 
$p\in M^2$.
If not, there exists a point $p_1\in M^2$
such that $\phi_v(p_1)\not \in C$
for all $v\in \R$.
Let $p_0$ be a point on $C$.
Since $M^2$ is connected,
we can take a continuous curve 
$p(s)$ ($0\le s\le 1$) 
in $M^2$ such that
$p(0)=p_0$ and $p(1)=p_1$.
Without loss of generality, we may assume that
$p(s)$ does not meet $C$ for $s>0$.
We set
$$
I_0:=\{s\in [0,1]\, ;\, 
\mbox{$\exists v\in \R$ such  that }\,\, \phi_{v}(p(s))\in  C\}. 
$$
Since $I_0$ contains $s=0$,
it is a
non-empty subset.
Since $\sigma'$ is transversal to 
$\xi$, $I_0$ is an open subset of $[0,1]$. 
Next, we show that $I_0$ is closed.
In fact, let $\{s_n\}_{n=1,2,3,...}$
be a sequence in $I_0$ converging to
a point $s_\infty\in (0,1]$.
Since $s_n \in I_0$,
there exists a point $q_n\in C$
and $v_n\in \R$ such that
$
\phi_{v_n}(q_n)=p(s_n).
$
Since $C$ is compact,
taking a subsequence if necessary,
we may assume that $\{q_n\}$ converges to $q_\infty\in C$.
Since $|df(\xi)|=1$, the map $v\mapsto f(\phi_{v}(q))$
gives a straight line parametrized by the arc-length parameter.
If $v_n\to 0$, we have $p(s_\infty)=q_\infty\in C$, in this case
$s_\infty\in I_0$, which contradicts the fact that
$p(s)\not \in I_0$ for $s>0$.
So we may assume that $v_n$ does not converge to $0$.
By replacing $\xi$ by $-\xi$ if necessary,
we may also assume that $v_n>0$ for all $n$.
Thus we have 
$
v_n=|f(q_n)-f(p(s_n))|.
$
So we have
$$
(v_\infty:=)\lim_{n\to \infty}v_n=|f(q_\infty)-f(p(s_\infty))|,
$$
and get the expression
$$
\phi^{}_{v_\infty}(q_\infty)=p(s_\infty),
$$
which implies $s_\infty\in C$.
Thus $I_0$ is open and closed in [0,1]
and so we can conclude that $I_0=[0,1]$,
which is a contradiction. Hence $C$
is a generator.

Next, we suppose that $f$ is not
$a$-orientable.
Then $\xi$ is not a globally defined
vector field on $M^2$, and
there exists a double covering 
$\pi:\hat M^2\to M^2$
such that $\xi$ can be considered as a
$C^\infty$-vector field of $\hat M^2$.
We set $C=\sigma(S^1)$.
We  now would like to show that 
the inverse image $\pi^{-1}(C)$
is connected.
If not, there are two disjoint closed curves $C_1$ and 
$C_2$ in $\hat M^2$ such that
$$
\pi^{-1}(C)=C_1\cup C_2.
$$
Then $C_i$ ($i=1,2$) is a closed regular curve on $M^2$
satisfying the condition of Lemma \ref{lem:transversal}.
Since $f\circ \pi$ is $a$-orientable,
$C_i$ ($i=1,2$)
must meet each asymptotic line.
We now fix an asymptotic line $L (\subset \hat M^2)$
as an integral curve of $\xi$.
Then there exist two points
$p,q\in \hat M^2$ such that
an arbitrarily fixed integral curve $\gamma$ of $\xi$
on $\hat M^2$ meets both of $\sigma_1$ and $\sigma_2$ at
$p$ and $q$ respectively.
Then $\pi(\gamma)$ meets $\sigma$ at 
$\pi(p)$ and $\pi(q)$.
Since $\sigma$ is a generator,
we have $\pi(p)=\pi(q)(=r)$.
Then the straight line $f\circ \pi(L)$ passes through 
the point $f(r)$ twice, which is a contradiction.
Thus $\pi^{-1}(C)$ is connected and
must be a generator of $f\circ \pi$.
Since $\pi^{-1}(C)$ meets each asymptotic line 
of $f\circ \pi$, $C$ also  
meets each asymptotic line 
of $f$, that is, $C$ is a generator.
\end{proof}

\begin{corollary}\label{cor:ori}
Let $f:M^2\to \R^3$ be an admissible developable 
p-frontal. Then the following three conditions are
mutually equivalent:
\begin{enumerate}
\item[{\rm (1)}] $f$ is $a$-orientable.
\item[{\rm (2)}] $M^2$ is orientable.
\item[{\rm (3)}] There exists a smooth unit asymptotic
vector field, which is globally defined on an generator.
\end{enumerate}
\end{corollary}

\begin{proof}
If $f$ is $a$-orientable, 
$M^2$ is orientable by Proposition~\ref{prop:dev_str}.
(Even when $f$ is non-co-orientable p-front,
the argument of the proof of Proposition~\ref{prop:dev_str}
works, since we do not use the existence of globally defined
unit normal vector $\nu$ there.)
Thus (1) implies (2).
Next we assume that $M^2$ is oriented. 
Let $\sigma(t)$ ($0\le t\le 1$) be a
generator of $f$.
We can choose $\xi(t)$ for each $t$
such that $(\sigma'(t),\xi(t))$ is positively oriented frame.
Then $\xi(0)=\xi(1)$ holds, which implies (3).
Finally, we assume (3).
Let $\sigma(t)$ ($0\le t\le 1$) be a
generator of $f$.
We can choose $\xi(t)$ for each $t$
such that $\xi(0)=\xi(1)$ holds.
Since $f$ is not $a$-orientable,
there exists a double covering
$\pi:\hat M^2\to M^2$ such that $f\circ \pi$ is $a$-orientable.
We denote by $C$ the image of $\sigma$.
Since $\xi(0)=\xi(1)$, $\pi^{-1}(C)$ consists of a union of
two closed curves $C_1$ and $C_2$.
However, we have seen in the last-half part of the proof of 
Lemma~\ref{lem:transversal}, it makes a
contradiction. Thus (3) implies (1).
\end{proof}

Here a function $\phi\in C^\infty(\R)$ satisfying $\phi(t+\pi)=-\phi(t)$
is called {\it anti-$\pi$-periodic}.
Firstly, we consider the  non-orientable case:

\begin{proposition}\label{prop:non-ori}
Let $f:S^1\times \R\to \R^3$ be a front given by
$$
f(t,v)=\hat\sigma(t)+v \hat\xi(t),\qquad
\left(\hat\sigma(t):=\int_0^t (a(s)\hat\xi(s)+b(s)\hat\eta(s))ds
\right),
$$
associated with a data $(a(t)dt,b(t)dt,\hat\xi(t),\hat \nu(t))$
on $S^1=\R/2\pi \Z$.
Then $f$ is
a double covering of a non-orientable
p-frontal $\bar f$
if and only if $\hat\xi(t),a(t)$ are anti-$\pi$-periodic.
Moreover 
\begin{enumerate}
\item[{\rm (1)}] $\bar f$ is co-orientable if and only if
$\nu(t)$ is $\pi$-periodic and $b(t)$ is anti-$\pi$-periodic.
\item[{\rm (2)}] $\bar f$ is non-co-orientable if and only if
$\nu(t)$ is anti-$\pi$-periodic and $b(t)$ is $\pi$-periodic.
\end{enumerate}
Conversely, an arbitrary non-orientable 
admissible developable p-frontal $\bar f$ has such a 
double covering.
\end{proposition}

\begin{proof}
The first assertion is obvious.
So we fix a non-orientable admissible developable p-frontal 
$\bar f:M^2\to \R^3$. 
Since $M^2$ is non-orientable, 
$f$ is not $a$-orientable.
Let $C$ be a generator of $\bar f$.
Then there exists a double covering 
$\pi:\hat M^2\to M^2$
such that $\xi$ can be considered as a
$C^\infty$-vector field on $\hat M^2$.
As seen in the proof of Lemma~\ref{lem:transversal},
$\pi^{-1}(C)$ is connected and
must be a generator of $\bar f\circ \pi$.
Let
$T:\hat M^2\to \hat M^2$ be the covering transformation.
We can take a parametrization $\gamma(t)$ 
($0\le t\le \pi$)
of $\pi^{-1}(C)$ such that
$\gamma(t+\pi)=T\circ \gamma(t)$.
If $f$ is co-orientable
(resp. non-co-orientable)
$\nu(t+\pi)=\nu(t)$ (resp. $\nu(t+\pi)=-\nu(t)$)
holds. In both of two cases, $\nu(t+2\pi)=\nu(t)$,
and so $\bar f\circ \pi$ is a front.
Thus, if we take the canonical coordinate 
(as in Remark \ref{rem:canonical})
with respect to $\gamma(t)$,
we get the assertion. 
\end{proof}

Like as in the non-orientable case, the following assertion holds.

\begin{proposition}\label{prop:non-co-ori}
Let $f:S^1\times \R\to \R^3$ be a front given by
$$
f(t,v)=\hat\sigma(t)+v \hat\xi(t),\qquad
\left(\hat\sigma(t):=\int_0^t (a(s)\hat\xi(s)+b(s)\hat\eta(s))ds
\right),
$$
associated with a data
$(a(t)dt,b(t)dt,\hat\xi(t),\hat \nu(t))$
on $S^1=\R/2\pi \Z$.
Then $f$ is
a double covering of a non-co-orientable
p-frontal $\bar f$
if and only if $\hat\nu(t)$ is anti-$\pi$-periodic.
Moreover 
\begin{enumerate}
\item[{\it (1)}] $\bar f$ is orientable if and only if
$\xi(t),a(t)$ is $\pi$-periodic and $b(t)$ is anti-$\pi$-periodic.
\item[{\it (2)}] $\bar f$ is non-orientable if and only if
$\xi(t),a(t)$ is anti-$\pi$-periodic and $b(t)$ is $\pi$-periodic.
\end{enumerate}
Conversely, an arbitrary non-co-orientable 
admissible developable p-frontal $\bar f$ has such a 
double covering.
\end{proposition}

\begin{proof}
The first assertion is obvious.
So we fix a non-co-orientable admissible developable p-frontal 
$\bar f:M^2\to \R^3$. 
Let $C$ be a generator of $\bar f$.
Then there exists a double covering 
$\pi:\hat M^2\to M^2$
such that $\nu$ can be considered as a
$C^\infty$-vector field on $\hat M^2$.
As seen in the proof of Lemma~\ref{lem:transversal},
$\pi^{-1}(C)$ is connected and
must be a generator of $\bar f\circ \pi$.
Let
$T:\hat M^2\to \hat M^2$ be the covering transformation.
We can take a parametrization $\gamma(t)$ 
($0\le t\le \pi$)
of $\pi^{-1}(C)$ such that
$\gamma(t+\pi)=T\circ \gamma(t)$.
If $f$ is orientable
(resp. non-orientable)
$\xi(t+\pi)=\xi(t)$ (resp. $\xi(t+\pi)=-\xi(t)$)
holds. In both of two cases, $\xi(t+2\pi)=\xi(t)$,
and so $\bar f\circ \pi$ is a-orientable.
Thus, if we take the canonical coordinate 
(as in Remark \ref{rem:canonical})
with respect to $\gamma(t)$,
we get the assertion. 
\end{proof}

\begin{proposition}\label{cor:umbilic}
There exists at least one asymptotic line
consisting of umbilic points on an 
admissible non-orientable 
developable p-front.
\end{proposition}

\begin{proof}
Firstly, we consider the case that 
$f$ is a front.
Since $M^2$ is non-orientable,
we can apply Proposition~\ref{prop:non-ori}.
Then $f$ can be represented by 
a quadruple
$(a(t)dt,b(t)dt, \hat \xi(t),\hat \nu(t))$
as in Corollary~\ref{cor:rep-1closed}, where 
$t\in S^1=\R/2\pi\Z$.
Suppose that the surface has no umbilics.
Then the Gaussian curve $\hat \nu$ 
is a regular curve by Proposition~\ref{prop:1-9}.
So we may assume that $\hat\nu=\hat\nu(t)$
is parametrized by the  arc-length parameter.
Then we have $\hat \xi=\pm \hat \nu'\times\hat \nu$.
Since $\pm$-ambiguity does not effect our construction,
we may set $\hat \xi=\hat \nu'\times\hat  \nu$.
Since 
$\hat \nu$ is $\pi$-periodic, 
so is  $\hat\nu'$.
Then we can conclude that
 $\hat \xi=\hat \nu'\times\hat  \nu$ is $\pi$-periodic,
which contradicts the anti-$\pi$-periodicity of $\hat \xi$.

Next, we consider the case that $f$ is a
non-co-orientable p-front.
By Proposition~\ref{prop:non-co-ori},
there exists a double covering
$
\pi:\hat M^2\to M^2
$
such that $f\circ \pi$ is
an admissible developable front, which is orientable.
Then 
$f$ can be represented by 
a quadruple
$(a(t)dt,b(t)dt, \hat \xi(t),\hat \nu(t))$
satisfying (2) 
of Proposition~\ref{prop:non-ori}.
Then the Gaussian curve $\hat \nu$ 
is a regular curve by Proposition~\ref{prop:1-9}.
So we may assume that $\hat\nu=\hat\nu(t)$
is parametrized by the  arc-length parameter.
We may set $\hat \xi=\hat \nu'\times\hat  \nu$.
Since $\hat \nu$ is anti-$\pi$-periodic, 
so is  $\hat\nu'$.
Then we can conclude that
 $\hat \xi=\hat \nu'\times\hat  \nu$ is $\pi$-periodic,
which contradicts the anti-$\pi$-periodicity of $\hat \xi$.
\end{proof}

\begin{remark}\label{rmk:cardioid}
For example, the cylinder over a cardioid 
$$
f(t,v)=(1-\sin t)\pmt{\cos t \\ \sin t \\ 0}+v \pmt{0\\0\\ 1}
$$
gives a developable front, which is non-co-orientable but 
orientable. The normal vector field of $f$ is given by
$$
\nu(t,v):=
\pmt{\left(\cos \frac{t}{2}+\sin \frac{t}{2}\right) (1-2 \sin t)
\\ \left(cos \frac{t}{2}-\sin \frac{t}{2}\right) (1+2\sin t)\\ 0}.
$$
In this case, the surface has no
umbilics. So in the above statement, 
non-orientability is crucial.
\end{remark}

Like as in Remark \ref{rem:infinite}, one can see
the existence of infinitely many
non-orientable p-fronts and non-co-orientable
p-fronts.
Finally, we remark that by modifying the argument
in the proof of Proposition~\ref{cor:umbilic}, one can prove that
{\it a developable M\"obius strip 
$($that is, a non-orientable developable immersion 
on $S^1\times (-1,1)\,)$
has at least one umbilic on it.}
(See \cite{KU}.)
It should be also remarked that
the fundamental properties of 
flat M\"obius strips are given in
\cite{H}, \cite{CK} and \cite{RR}. 
Recently, the existence of real analytic
flat M\"obius strips of a given isotopy type,
whose centerline is
a geodesic or a line of curvature, was shown in
\cite{KU}.
Also, in \cite{KU}, examples of weakly complete
flat fronts containing a M\"obius strip
were given. These fronts cannot be complete,
since complete flat fronts are orientable
 (cf. Theorem~B).

\section{Completeness}
In the introduction, we have defined two types of
completeness for
flat fronts. In this section, we 
determine the structure
of complete flat fronts.
Firstly, we prepare two lemmas.

\begin{lemma}\label{lem:3-0}
A complete flat front $f:M^2\to \R^3$ 
is weakly complete.
\end{lemma}

\begin{proof}
Let $ds^2=df\cdot df$ be the first fundamental form.
Since $f$ is complete, there exists a symmetric tensor
$T$ having compact support $K:=\op{supp}(T)$ 
such that $g:=ds^2+T$ is 
a complete Riemannian metric on $M^2$.
It is well-known that 
any two Riemannian metrics $g_1,g_2$
on $M^2$ satisfies $cg_2<g_1< (1/c)g_2$
on a given compact set for suitable positive number $c$.
So there exists
a constant $0<c<1$ such that
$$
ds^2_\#:=df\cdot df+d\nu\cdot d\nu>cg\quad \mbox{on $K$}.
$$
On the other hand, it is obvious that
$ds^2_\#\ge ds^2$ holds on $M^2\setminus K$.
So $ds^2_\#>cg$ holds on $M^2$.
Since $g$ is complete, so is $ds^2_\#$. 
\end{proof}

\begin{lemma}\label{lem:3-1}
A weakly complete flat front $f:M^2\to \R^3$ which has
no umbilics is an $a$-complete developable front.
\end{lemma}

\begin{proof}
By Proposition~\ref{ftodev}, $f$ is a
developable front. 
So it is sufficient to show that
each asymptotic line is complete.
We take an asymptotic line $\gamma(t)$ ($t\in \R$) arbitrarily.
Then $\nu$ is constant along  $\gamma(t)$, that is,
$d\nu=0$ on $\gamma$, and
$ds^2_\#$-length of $\gamma$ is equal to that of
$ds^2$-length.
Since $ds^2_\#$-length is greater than or equal to 
$ds^2$-length in general, we can conclude that $\gamma$
is a geodesic of $(M^2,ds^2_\#)$.
Since $ds^2_\#$ is a complete metric, 
$\sigma$ has infinite length and the image of $f\circ \gamma$
must coincide with a complete straight line.
\end{proof}

The authors do not know whether weakly complete flat fronts
are developable (more strictly, admissible) or not.
Later, we will show that complete flat fronts are all 
admissible developable. (See the proof of Theorem~A in this
section.)
Conversely, we can prove the following
(see also Corollary~\ref{4-2}):

\begin{proposition}\label{prop:3-2}
Let $f:M^2\to \R^3$ be an admissible $a$-complete
developable front which has a closed generator.
Then $f$ is weakly complete. 
\end{proposition}

\begin{proof}
Taking the double covering, we may assume 
$f$ is $a$-orientable. (See Section 3.)
By Corollary~\ref{cor:rep-1closed}, 
$M^2$ can be identified with $S^1\times \R$,
and
$f$ is constructed from
a quadruple $(a(t)dt,b(t)dt,\hat \xi(t),\hat\nu(t))$
on $S^1=\R/2\pi\Z$.
Then the
lift metric is given by (cf. \eqref{eq:lift})
$$
ds^2_\#=df\cdot df+d\nu\cdot d\nu
=\left(a^2+(b+v \mu_1)^2+(\mu_2)^2\right)dt^2+2a dt dv+ dv^2,
$$
where 
$
\hat\xi'(t)=\mu_1(t)\hat \eta(t)$ and $
\hat\nu'(t)=\mu_2(t)\hat \eta(t)$. 
Since the symmetric matrix
$$
\pmt{a^2+(b+v \mu_1)^2+(\mu_2)^2 & a \\ a & 1}
$$
has positive eigenvalues $k_1,k_2$
($k_2>k_1$) such that
$$
k_1+k_2=1+a^2+(b+v \mu_1)^2+(\mu_2)^2,\quad 
k_1k_2=(b+v \mu_1)^2+(\mu_2)^2,
$$
we have
$$
k_2\ge k_1=\frac{(b+v \mu_1)^2+(\mu_2)^2}{k_2}
\ge \frac{(b+v \mu_1)^2
+(\mu_2)^2}{1+a^2+(b+v \mu_1)^2+(\mu_2)^2}.
$$
We set $A(t,v):=(b(t)+v \mu_1(t))^2+\mu_2(t)^2$.
If $A$ has a zero at $(t_0,v_0)$,
then we have $\mu_2(t_0)=0$ and $b(t_0)+v \mu_1(t_0)=0$.
In particular, $\hat\nu'(t_0)=\mu_2(t_0)\eta(t_0)=0$.
Since $f$ is a front and $\hat\nu'(t_0)=0$, 
\eqref{front-cond} holds by
Theorem~\ref{thm:w-rep}.
Then  we have that
$$
0\ne b(t_0)=-v \mu_1(t_0)=0,
$$
which makes a contradiction.
Thus $A(t,v)>0$ on $tv$-plane. Since $\dy\lim_{v\to \infty} A(t,v)$
diverges for all $t\in S^1$, the positive number
$$
m:=\min_{(t,v)\in \R^2}A(t,v)
$$
is well-defined.
Since $x\mapsto x/(1+a^2+x)$ is monotone increasing,
we have $k_1\ge m/(1+a^2+m)$.
Then we have that
$$
ds^2_\# \ge \frac{m}{1+a^2+m} (dt^2+dv^2).
$$
Since $dt^2+dv^2$ gives the complete metric on $S^1\times \R$,
so does $ds^2_\#$, which proves the assertion.
\end{proof}

\begin{proposition}\label{prop:no-umbilics}
$f:M^2\to \R^3$ be an $a$-complete and $a$-orientable
developable front whose singular set is 
non-empty and compact.
Suppose that there are no umbilics on $M^2$.
Then $f$ is admissible and the singular set 
of $f$ is a generator.
\end{proposition}

\begin{proof}
There exists a unit asymptotic vector field 
$\xi$ on $M^2$.
Since the umbilic set is empty,
the principal radius function $\rho$ 
is a real-valued $C^\infty$-function on $M^2$ (see Section 1).
Thus we can define a $C^\infty$-function by
$
\phi=d\rho(\xi):M^2\to \R.
$
We fix a point $p\in M^2$ arbitrarily.
We take a curvature line 
coordinate neighborhood
$(U;u,v)$ ($|u|<\epsilon$,\,\,$v\in \R$) centered at $p$
satisfying properties (1)--(2) given in Section 1.
Then by Lemma \ref{lem:gmassey}, there exist 
$C^\infty$-functions $A(u), B(u)$ such that
\begin{equation}\label{a1}
\rho(u,v)=A(u)v+B(u).
\end{equation}
By definition, we have
\begin{equation}\label{a2}
\phi(u,v)=A(u).
\end{equation}
Let $Z_{\phi}$ be the zeros of $\phi$ on $M^2$.
Suppose that $Z_{\phi}$ is non-empty and is not equal to
$M^2$. Then we may assume $p=(0,0)$ is
on the boundary of $Z_\phi$.
Then there exists a sequence $\{(u_n,v_n)\}_{n=1,2,3,...}$ in $U$
such that $\dy\lim_{n\to \infty} (u_n,v_n)=(0,0)$ and 
$
A(u_n)=\phi(u_n,v_n)\ne 0.
$
Then we have
$$
\lim_{n\to \infty}A(u_n)=\lim_{n\to \infty}\phi(u_n,v_n)=\phi(p)=0.
$$
Since $A(u_n)\ne 0$, we can set
$
w_n:=-{B(u_n)}/{A(u_n)}.
$
Then $\rho(u_n,w_n)=0$ holds, that is,
 $\{(u_n,w_n)\}$ is a sequence consisting of singular points
in $U$. If $B(0)\ne 0$, then this singular set is unbounded,
which contradicts the compactness of the singular set.
So we can conclude that $B(0)=0$. Then
$
\rho(p)=A(0)v+B(0)=0,
$
which implies that the curve $v \mapsto (0,v)$
consists of singular points, which
contradicts the compactness of the singular set again.
Thus we have $Z_{\phi}=M^2$, since
we have assumed that $Z_{\phi}$ is non-empty.
However, if $Z_{\phi}=M^2$, then
$
\rho(u,v)=B(u),
$
that is, $\rho$ does not depend on $v$.
If $\rho$ has zeros, then the singular set contains an asymptotic line.
This contradicts the compactness of the singular set.
So $\rho$ never vanishes. 
Then this implies that there are no singular points on 
$M^2$, which contradicts our assumption. 
Thus $Z_{\phi}$ must be an empty set. 
In particular, $A(u)$ never vanishes,
and so each asymptotic line has a  singular point
exactly at $v=-B(u)/A(u)$.

Now, we fix a point  $p\in M^2$ arbitrarily.
We have an expression \eqref{a1}.
Since we have seen that $Z_{\phi}$ is the empty set,
$A(u)$ never vanishes, and so
the singular set $S(f)$ on $(U;u,v)$ can be parametrized by
$
v:=-{B(u)}/{A(u)},
$
which implies that the singular curve is transversal
to the asymptotic direction (i.e. $v$-direction)
and the singular point is unique in each asymptotic line.
Since we have already seen that each asymptotic line 
contains a singular point, we can conclude that $S(f)$
is a generator.
\end{proof}

\begin{corollary}\label{cor:xig}
Let $f$ be as in Proposition~\ref{prop:no-umbilics}.
Then the Gaussian curve and the asymptotic Gauss map
are both regular curves without inflection points.
\end{corollary}
\begin{proof}
We may assume that $M^2$ is connected,
and
can take the singular set $S(f)$ as
a generator.
Then $M^2$ is diffeomorphic to $S^1\times \R$ by 
Proposition~\ref{prop:dev_str}, and we can take
a quadruple $(a(t)dt,b(t)dt,\hat \xi,\hat \nu)$
which represents $f$.
Here $b(t)$ vanishes identically.
In fact, if $b(t_0)\ne 0$, 
\begin{equation}\label{eq:b}
f_t(t,v)\times f_v(t,v)
=(b(t)+v \mu_1(t))\hat\eta(t)\times \hat\xi(t)
\end{equation}
does not vanish at $(t_0,0)$, which contradicts the
fact that $t\mapsto (t,0)$ is the singular curve.
Since $f$ has no umbilics,
$\hat\nu(t)$ is a regular curve.
Suppose that $\hat\xi'(t_0)=0$, that is, 
$\mu_1(t_0)=0$.
Then \eqref{eq:b} implies that
$v\mapsto (t_0,v)$ is a singular curve,
which contradicts 
the compactness of the singular set.
So $\hat \xi$ must be a regular curve.
\end{proof}

\begin{theorem}\label{4-4}
Let $M^2$ be a connected 2-manifold, and
$f:M^2\to \R^3$ an orientable
weakly complete flat front whose singular set is 
non-empty and compact.
Then $f$ has no umbilics.
\end{theorem}

\begin{proof}
We suppose that $\mc U_f$ is non-empty.
Since any asymptotic line passing through a non-umbilic
point never meets the 
umbilic set $\mc U_f$ (cf. Lemma \ref{lem:gmassey}),
Lemma \ref{lem:3-1}
yields that  
the restriction $\hat f:=f|_{M^2\setminus \mc U_f}$
is an $a$-complete developable front.
Since $\mc U_f$ does not contain any singular points,
$\hat f$ is $a$-complete having compact
singular set.
If $M^2=\mc U_f$, then $f$ has no singular points, 
a contradiction.

So there exists a boundary point $p$ on $\mc U_f$.
Then there exists a sequence $\{q_n\}$ on $M^2\setminus \mc U_f$
converging to $p$.
Let $L_n$ be the asymptotic line on $M^2$ passing through $q_n$.
By Proposition~\ref{prop:no-umbilics} for $\hat f$,  
the singular set is a generator, and
the condition (3) of Definition~\ref{def:sigma} yields that 
there exists a unique singular point $s_n$ on $L_n$.
Since the singular set is compact, we may assume that $s_n$
converges to a point $s$.
Since $s$ is a singular point, $s\not\in \mc U_f$,
and there exists an asymptotic line $L$ passing through $s$.
Since $L_n$ must converge to $L$ and $q_n\in L_n$, 
the asymptotic line $L$ meets an umbilic point $p$, 
 which contradicts Lemma \ref{lem:gmassey} again. 
Thus $\mc U_f$ must be the empty set.
\end{proof}

\medskip
\noindent
{\it Proof of Theorem~A }:
Since $f$ is complete, it is weakly
complete by Lemma \ref{lem:3-0}.
First, we consider the case that
$f$ is an orientable front.
By Theorem~\ref{4-4},
$f$ has no umbilics.
Then by Theorem~\ref{prop:no-umbilics},
$f$ is admissible and the singular set is a generator.
Since $M^2$ is connected,
$M^2$ is diffeomorphic to $S^1\times \R$ by 
Proposition~\ref{prop:dev_str}.

Next we consider the case that $f$ is a non-orientable
front or a non-orientable p-front. By 
Proposition~\ref{prop:non-ori},
we can take a double covering
$\pi:\hat M^2\to M^2$ such that
$\hat M^2$ is orientable and $f\circ \pi$
is a complete flat front.
Then, by Theorem~\ref{4-4},
$f\circ \pi$ (and also $f$) has no-umbilics.
On the other hand, Proposition~\ref{cor:umbilic} yields that
$f$ (and so $f\circ \pi$ also) 
has an umbilical point, which makes a contradiction.  
So this case is impossible.

Finally we consider the remaining case 
that $f$ is a non-co-orientable
(but orientable) p-front. 
Then we can take a double covering
$\pi:\hat M^2\to M^2$ such that
$f\circ \pi$ is a  complete flat front.
Since $M^2$ is orientable, so is $\hat M^2$.
Then $f$ is $a$-orientable by Corollary~\ref{cor:ori}. 
Hence, by Theorem~\ref{4-4},
$f\circ \pi$ has no umbilics, and
is an admissible developable front
such that the singular set is a generator of $f\circ \pi$. 
Thus by Corollary~\ref{cor:rep-1closed}, 
$\hat M^2$ is diffeomorphic to $S^1\times \R$.
Moreover, since the singular curve
is a generator, 
the proof of Corollary \ref{cor:xig} yields that
$f\circ \pi$ 
can be represented by a quadruple $(a(t)dt,0,\hat\xi,\hat\nu)$ 
such that $a(t),\hat\xi(t)$ are $\pi$-periodic and 
$\hat\nu(t)$ is anti-$\pi$-periodic, 
that is, $\hat\nu(t+\pi)=-\hat\nu(t)$ (cf. 
Proposition~\ref{prop:non-co-ori}). 
Then $f$ satisfies
\begin{equation}\label{eq:f-period}
f(t,v)=f(t+\pi,v) \qquad (t,v\in \R).
\end{equation}
In particular,
the curve 
$$
\sigma_1:S^1\ni t\mapsto (t,1)\in S^1\times \R(=\hat M^2)
$$
is $\pi$-periodic.
Since $b(t)$ vanishes identically and
$\mu_1(t)$ never vanishes by Corollary~\ref{cor:xig}, 
the equation $v=0$ characterizes the singular set 
(cf. Theorem~\ref{thm:w-rep}).
Thus $\sigma_1(t)$
does not meet the singular set of $f$, and
$\hat \sigma_1(t)=f\circ \sigma_1(t)$ is a regular
space curve.
Then $\hat\nu(t)$ or $-\hat\nu(t)$ must coincide with
$$
{\hat\xi(t)\times \hat\sigma'_1(t)}
/{|\hat\xi(t)\times \hat\sigma'_1(t)|}.
$$
Since $\hat\xi(t),\hat\sigma_1(t)$ are both $\pi$-periodic, so is
$\hat\nu(t)$, which contradicts the anti-$\pi$-periodicity of
$\hat\nu(t)$. So this case also never happens.
\qed

\begin{corollary}\label{4-2}
$f:M^2\to \R^3$ be a weakly complete
front whose singular set is non-empty and compact.
Then $f$ is complete.
\end{corollary}

\begin{proof}
Without loss of generality, $M^2$ is connected.
By Theorem~\ref{4-4}, $f$ has no umbilics.
Then by Theorem~\ref{prop:no-umbilics},
$f$ is admissible and the singular set is a generator.
Since $M^2$ is connected,
Corollary~\ref{cor:rep-1closed} yields that
$f$ can be represented by a
quadruple $(a(t)dt,0,\hat\xi,\hat \nu)$ ($t\in S^1$).
Thus the first fundamental form $ds^2$
of $f$ can be given by
$$
ds^2=\left(a(t)^2+v^2\mu_1(t)^2\right)dt^2+2a(t) dt dv+dv^2,
$$
where $\hat \xi'=\mu_1\hat \eta$
\,\,($\hat \eta=\hat \xi\times \hat \nu$).
Let $\phi:\R\to [0,1]$ be a smooth function
such that $\phi(v)=1$ for $v\in [-1,1]$
and $\phi(v)=0$ for $|v|\ge 2$. 
It is sufficient to show that
the metric
$$
g:=ds^2+\phi(v)^2\mu_1(t)^2dt^2
$$
is a complete Riemannian metric on $S^1\times \R$.
Since the singular set is compact and the zeros of
$\mu_1$ correspond to the linear singularity (cf. 
Definition~\ref{def:linear}), 
$\mu_1(t)$ never vanishes on $S^1$,
and we may set
$$
m:=\min_{t\in S^1}|\mu_1(t)|>0.
$$
Then we have
\begin{equation}\label{eq:metric}
g\ge (a^2+m^2)dt^2+2a dt dv+dv^2.
\end{equation}
Since the two eigenvalues $k_1,k_2$ 
($k_2\ge k_1$) of
the matrix
$
\pmt{a^2+m^2 & a \\ a & 1}
$
satisfy $k_1k_2=m^2$
and $k_2\le k_1+k_2=a^2+m^2$,
we have 
$$
\min_{t\in S^1}k_1(t)=
\frac{m^2}{\dy\max_{t\in S^1}k_2(t)}
\ge \frac{m^2}{A^2+m^2}>0,
$$
where $A:=\max_{t\in S^1}|a(t)|$.
In particular, the right hand side 
of \eqref{eq:metric}
gives a complete Riemannian metric on the $tv$-plane. 
Thus $g$ is also complete.
\end{proof}

\medskip
\noindent
{\it Proof of Theorem~B}:
Let $M^2$ be a connected 2-manifold, and
$f:M^2\to \R^3$  a complete flat front.
Then as seen in the proof of Theorem~A, $M^2$ is
orientable and
$f$ is admissible, where the singular set  
is a generator.
Then $f$ can be represented by
a quadruple $(a(t)dt,b(t)dt,\hat \xi,\hat \nu)$
($t\in S^1$)
associated with the singular set as a generator.
Then as seen in the proof of Corollary~\ref{cor:xig},
$b(t)$ vanishes identically.
Thus,  \eqref{eq:rep-2} reduces to the formula in Theorem~B.

Conversely, let $f$ be a map given by
the formula in Theorem~B.
Since $\hat \xi(t)$ is a regular curve without inflection
points, so is the dual spherical curve $\hat \nu(t)$.
(The singular point of $\hat\xi(t)$ (resp.~$\hat\nu(t)$)
corresponds to the inflection point of $\hat \nu(t)$ 
(resp. $\hat\xi(t)$)).
In particular  $\nu'(t)$ does not vanishes,
namely $f$ is a front without umbilical points.  
Moreover, $f$ is weakly complete by 
Proposition~\ref{prop:3-2}.
By the construction of $f$, 
its singular set is compact,
and $f$ is complete by Corollary~\ref{4-2}.
\qed

\begin{example} \label{cone}(A cone over a locally convex
spherical curve)
Let $\hat \xi:S^1\to S^2$ be 
a spherical closed regular  curve 
which has no inflection points.
By setting $a(t)=0$ in Theorem~B,
the map
$$
f(t,v):=v\hat \xi(t)
$$
gives a complete flat front.
The image of the singular curve $\{v=0\}$ in the
 $(t,v)$-plane
consists of an isolated  cone-like singular point.
We remark that
non-zero principal curvature lines 
on $f$ are all periodic
(cf. Proposition~\ref{prop:p-period}).
\end{example}

\begin{definition}
Let $J$ be an open interval and
$\gamma:J\to S^2$ a regular curve.
Then $\gamma(t_0)$ is called a {\it vertex} if
$t=t_0$ is a critical point of the geodesic curvature 
function $\kappa_g(t)$ of $\gamma(t)$.
A vertex $t=t_0$ is called {\it generic} if $\kappa''_g(t_0)\ne 0$.
\end{definition}
As an application of Theorem~B and 
Proposition~\ref{prop:singular},
we prove the following:
\begin{corollary}\label{cor:sing}
Let $f:S^1\times \R\to \R^3$ be a 
complete flat front whose singular 
set is non-empty,
and $\hat \nu$ its Gaussian curve.
Let $\Gamma$ be the set of
real numbers $\delta$
so that $f_\delta$ has a singular point
which is not a cuspidal edge nor a swallowtail.
\begin{itemize}
\item[(a)] If $f$ is real analytic, then $\Gamma$ is
discrete.
\item[(b)] If $\hat \nu$ admits only generic vertices,
then $\Gamma$ is bounded.
\end{itemize}
In particular,
$\Gamma$ is finite if $f$ satisfies both $($a$)$ and $($b$)$.
\end{corollary}

\begin{proof}
Since $f$ has a singular point, we may assume that
$f$ is associated with the quadruple 
$(a(t)dt,0,\hat \xi,\hat \nu)$.
Since $\hat \nu(t)$ is a regular curve
(cf. Corollary \ref{cor:xig}), we may assume that
$t$ is the arc-length parameter. 
Now, we use the notations as in \eqref{eq:mus}.
Then $\mu_2(t)=1$ identically, and
$-\mu_1(t)$ is just the non-vanishing
geodesic curvature function of
the spherical curve
$\hat \nu(t)$, since
$\hat\nu'=\hat\eta$.
We set
$$
\phi_\delta(t):=a(t)- \delta \left(\frac{1}{\mu_1(t)}\right)'
\qquad (\prime=d/dt).
$$
Since $(a(t)dt,\delta dt,\hat \xi,\hat \nu)$
represents $f_\delta$ (cf.
Proposition~\ref{prop:singular}),
a singular point $(t,v)\in S^1\times \R$ of $f_\delta$ 
which is neither a cuspidal edge nor
a swallowtail appears only when 
$\phi_\delta(t)=\phi'_\delta(t)=0$.
In this case, $(t,-\delta/\mu_1(t))$ gives such a
singular point.

We shall now prove (a).
Suppose that $f$ is real analytic and
there exists a bounded
sequence $\{\delta_n\}_{n=1,2,3,...}$ 
consisting of distinct points
such that 
$f_{\delta_n}$ has a singular point 
$(t_n,v_n)\in S^1\times \R$ 
which is neither a cuspidal edge nor
a swallowtail.
Then we have $v_n=-\delta_n/\mu_2(t_n)$.
Since $t_n\in S^1$, we may assume that $\{t_n\}_{n=1,2,3,..}$ 
converges to a point $t_\infty\in [0,2\pi)$.
We set 
$$
r(t):=1/\mu_1(t)>0.
$$
Since $\phi_{\delta_n}(t_n)=\phi'_{\delta_n}(t_n)=0$,
the function 
$$
\psi(t):=r'(t)a'(t)-r''(t)a(t)
$$
vanishes at $t=t_n\,\, (n=1,2,3,...)$.
Since $\psi(t)$ is real analytic, this implies that
there exists a positive integer $N$ such that
$t_n=t_N$ for $n\ge N$.
Then, for $n\ge N$ we have
$$
0=\phi_{\delta_n}(t_n)=\phi_{\delta_n}(t_N)=
a(t_N)- \delta_n 
\left. \left(\frac{1}{\mu_1(t_N)}\right)'\right|_{t=t_N},
$$
which contradicts that $\{\delta_n\}$
consists of mutually distinct points. 
This proves (a).
Next, suppose now that $\hat \nu$ has only generic vertices.
Then $|r'(t)|+|r''(t)|$ never vanishes for all $t\in S^1$,
since $-1/r(t)$ is the geodesic curvature of $\hat \nu(t)$.
If $f_\delta$ admits a singular point $(t_0,v_0)$
which is neither a cuspidal edge nor
a swallowtail, we have 
$a(t_0)=\delta r'(t_0)$ and $a'(t_0)=\delta r''(t_0)$.
In particular, 
$$
|\delta|=
\frac{|a(t_0)|+|a'(t_0)|}
{|r'(t_0)|+|r''(t_0)|}
\le 
\frac{\max_{t\in S^1}(|a(t)|+|a'(t)|)}
{\min_{t\in S^1}(|r'(t)|+|r''(t)|)},
$$
which implies that $\Gamma$ is bounded.
\end{proof}
\begin{remark}
A similar result for
flat fronts in $H^3$ 
is given in \cite{KRSUY}.
\end{remark}

\section{Embeddedness of ends}
In this section, we shall prove Theorem~C
and Theorem~D in the introduction.
Let $f$ be a complete flat front associated with 
a pair $(a(t)dt,\hat\xi(t))$ on $S^1$
as in Theorem~B.
Without loss of generality, we may assume that
$t$ is the arclength parameter of $\xi(t)$ (cf.~Corollary \ref{cor:xig}).
Then $\xi(t)$ and $\xi'(t)$ are linearly independent.
First, we suppose $f(t,v)$ ($t\in S^1=\R/2\pi\Z,\,\, v\in \R$)
has an end with self-intersection.
There are two ends $\{v>0\}$ and $\{v<0\}$.
Without loss of generality, we may 
assume that the end having self-intersections
is $\{v>0\}$.
Then there exist $(u_n,v_n),(x_n,y_n)\in S^1\times \R$ 
such that
\begin{equation}\label{intersection}
f(u_n,v_n)=f(x_n,y_n)\qquad (x_n\ne u_n,\,\, 0<v_n<y_n),
\end{equation}
and the sequence $\{v_n\}$ diverges to $\infty$. 
Taking a subsequence if necessary, we may assume that 
$u_n\to u_\infty$ and $x_n\to x_\infty$, where
$u_\infty,x_\infty\in [0,2\pi)$.
Using the expression \eqref{eq:B}, the singular set
is parametrized by a curve
$$
\hat\sigma(t):=\int_0^t \hat\xi \alpha=f(t,v)-v \hat\xi(t).
$$
Since $|\hat\xi(u_n)|=|\hat\xi(x_n)|=1$ and
$|\hat\xi(u_n)\cdot\hat\xi(x_n)|\le 1$,
the function
$v\mapsto |\hat\xi(u_n)-v \hat\xi(x_n)|^2$
is a monotone increasing function for $v\ge 1$.
Since  $y_n/v_n\ge 1$,
we have by \eqref{intersection}, 
$$
\left|\frac{\hat\sigma(x_n)-\hat\sigma(u_n)}{v_n}\right |
=\left|\hat\xi(u_n)-\frac{y_n}{v_n} \hat\xi(x_n)\right|
\ge |\hat \xi(u_n)-\hat \xi(x_n)|.
$$
Since the singular curve $\hat \sigma$ is bounded,
if we let $n\to \infty$, then $v_n\to \infty$
and the right hand side converges to $0$.
Thus we have
$
\hat\xi(u_\infty)=\hat\xi(x_\infty).
$

On the other hand, it holds that
\begin{align*}
\frac{1}{v_n}\left(\frac{\hat\sigma(x_n)-\hat\sigma(u_n)}{x_n-u_n}
\right)
+\frac{\hat\xi(x_n)-\hat\xi(u_n)}{x_n-u_n}
&=\frac{f(x_n,v_n)-f(x_n,y_n)}{v_n(x_n-u_n)} \\
&=\left(\frac{1-y_n/v_n}{x_n-u_n}\right)\hat\xi(x_n).
\end{align*}
Suppose now that $u_\infty=x_\infty$.
Since the singular curve $\hat \sigma$ is bounded,
the left hand side converges to $\hat\xi'(u_\infty)$ when $n\to \infty$,
and the right hand side also must converge.
In particular, $\frac{1-y_n/v_n}{x_n-u_n}$ must converge,
and $\hat\xi'(u_\infty)$ must be proportional to 
$\hat\xi(u_\infty)$, which contradicts
the fact that $\hat\xi(u_\infty),\hat\xi'(u_\infty)$ are linearly
independent.
Thus we have $u_\infty\ne x_\infty$, and so
$\hat\xi$ must have a self-intersection.
By the duality,  the Gaussian curve $\hat\nu$ 
also has a self-intersection,
since the dual of a convex spherical curve is also convex.
(Here we used the fact that 
$\hat\xi$ and $\hat \nu$ have no inflections).

Next, we consider 
the intersection of the front $f$ with
the sphere $S^2(r)$ of (sufficiently large) 
radius $r$ centered at
the origin.
Rescaling the surface $f$, 
it is equivalent to 
consider the intersection 
between the unit sphere $S^2(1)$
and the flat front 
$$
F_\delta(t,v):=\delta \hat \sigma(t)+v \hat \xi(t)
\qquad (t\in S^1,\,\, v\in \R),
$$
corresponding to the data $(\delta\alpha,\hat \xi)$, 
where $\delta=1/r$.
Let $\Gamma_\delta$ be the intersection of
the image of $F_\delta$
with the unit sphere $S^2(1)$.

Let $k$ be a constant so that $0<k<1/\sqrt{2}$.
Now we assume
$$
|\delta|<k/m \qquad (m:=\max_{t\in S^1}|\hat\sigma(t)|).
$$
Since 
$$
|\hat\sigma\cdot \hat \xi|\le m,\qquad 
\left||\hat\sigma|^2-(\hat\sigma\cdot \hat \xi)^2
\right|
=|\hat\sigma\times \xi|\le m,
$$
we have that
\begin{align*}
&-\delta\hat\sigma\cdot \hat \xi
+ \sqrt{1-\delta^2
\left(|\hat\sigma|^2-(\hat\sigma\cdot \hat \xi)^2\right)}
>-k+\sqrt{1-k^2}(>0), \\
&-\delta\hat\sigma\cdot \hat \xi
- \sqrt{1-\delta^2
\left(|\hat\sigma|^2-(\hat\sigma\cdot \hat \xi)^2\right)}
<k-\sqrt{1-k^2}(<0).
\end{align*}
Then we have
\begin{align*}
C_{k,\delta}&:={F_{\delta}}^{-1}(\Gamma_\delta)\cap \{v>0\}\\
&=\left\{(t,v)\,;\, 
v=-\delta\hat\sigma(t)\cdot \hat \xi(t)
+ \sqrt{1-\delta^2
\left(|\hat\sigma(t)|^2-(\hat\sigma(t)\cdot \hat \xi(t))^2\right)}
\right\}.
\end{align*}
This expression tells us that 
$C_{k,\delta}$
is a simple closed regular curve on $S^1\times \R$.
(This proves also that a 
complete flat front is a proper mapping.)
Now, we suppose that $f$ has an embedded end
$\{v>0\}$.
Then for sufficiently large fixed $k$,
the image $F_{\delta}(C_{k,\delta})$ ($0<\delta<k$) 
is a simple closed 
regular curve on $S^2(1)$.
By taking the limit $\delta\to 0$ in this expression,
$C_{k,\delta}$ converges to a curve $\{v=1\}$,
and we have 
$$
\lim_{\delta\to 0}F_{\delta}(C_{k,\delta})
=\hat \xi(S^1).
$$
Here, the above expression of $C_{k,\delta}$ and its image 
$F_{\delta}(C_{k,\delta})$ are both meaningful also for $\delta\le 0$,
and depend on $\delta$ smoothly at $\delta=0$.
Since $\hat \xi$ has non-vanishing geodesic curvature,
so does $F_{\delta}(C_{k,\delta})$ for sufficiently small
$\delta$. 
Since it is well-known that
a simple closed spherical curve having
non-vanishing geodesic curvature 
is a convex curve on an open hemi-sphere,
we can conclude that $F_{\delta}(C_{k,\delta})$ is
a convex curve on the unit sphere.
Thus the limit curve
$\hat \xi(S^1)$ must be also a convex curve.
Since $\hat \nu$ is the dual of $\hat \xi$,
$\hat \nu$ is also a convex curve.
\qed

\begin{figure}[htb]
  \begin{center}
         \includegraphics[width=5cm]{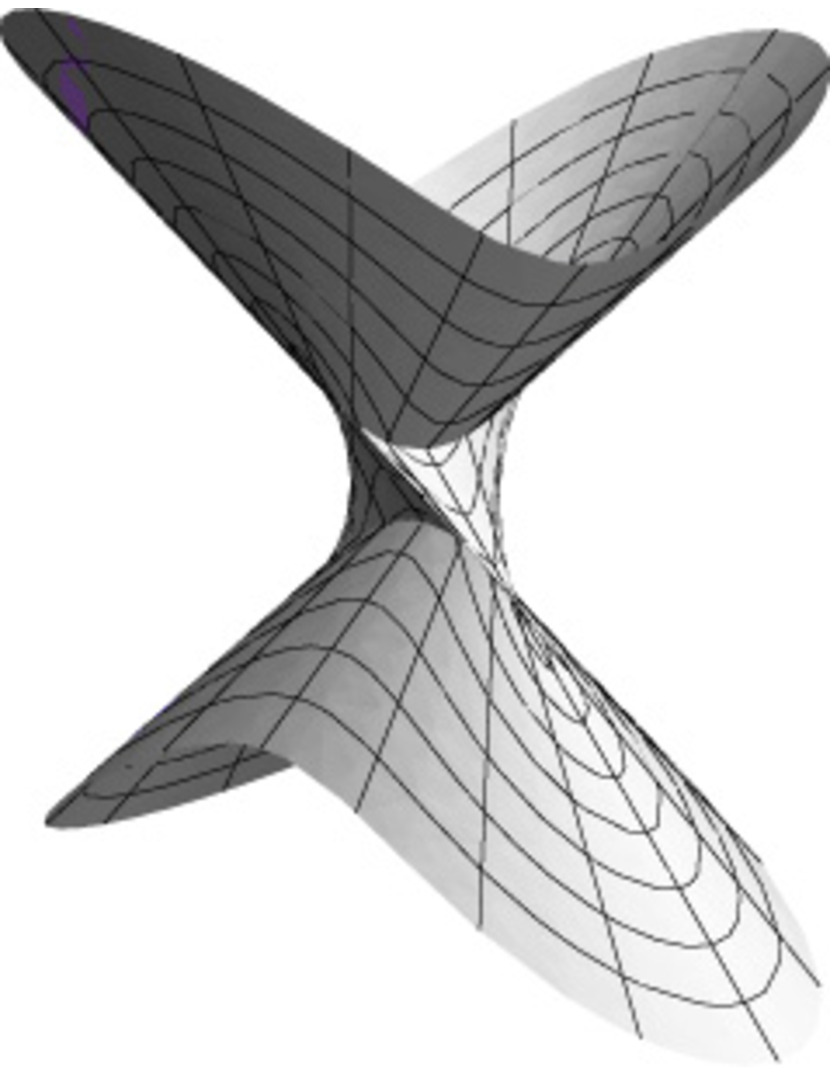}\qquad 
         \includegraphics[width=4cm]{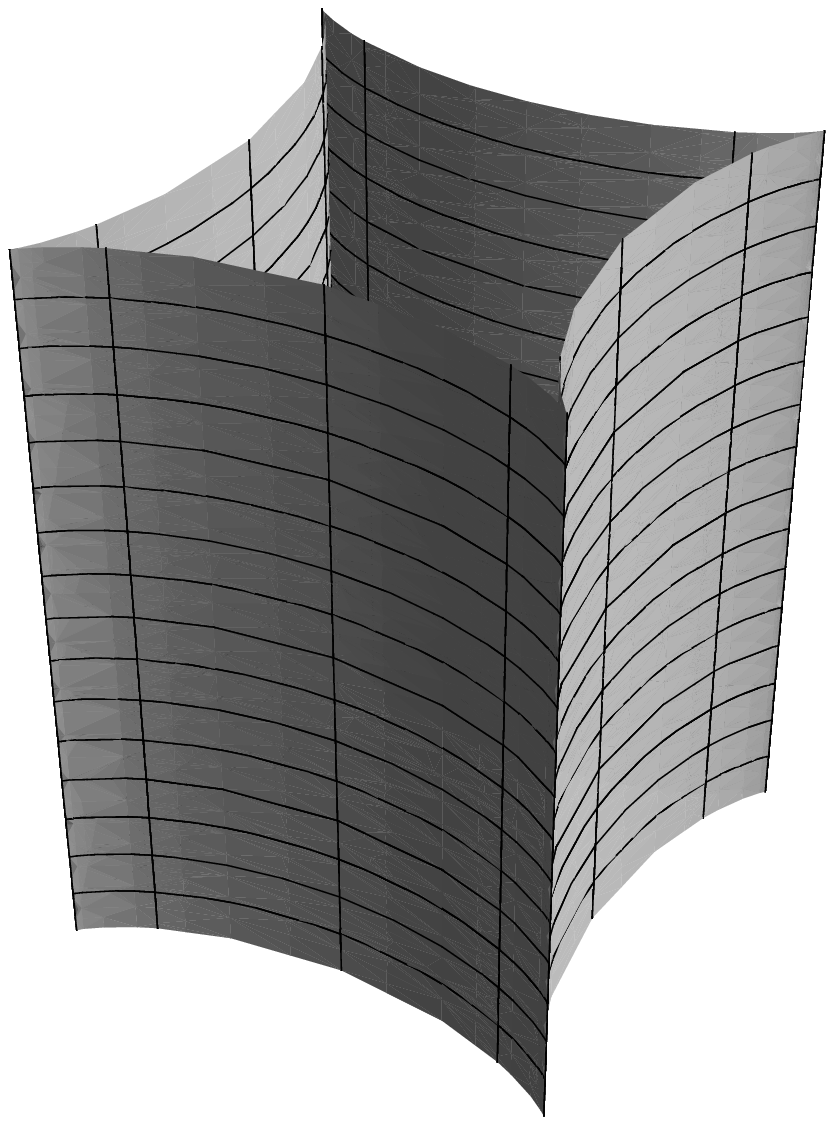}
\caption{The flat surface for $n=2,\,\, \phi=\pi/4$ and its caustic}
\label{hamen}
  \end{center}
\end{figure}

Finally, we give a non-trivial example with embedded ends.
 
\begin{example}
Let $\phi$ be a real number such that $|\phi|<\pi/2$.
We set
$$
\xi_{\phi}(t)=(\cos t \cos \phi,\sin t \cos \phi,\sin \phi)
\qquad (0\le t\le 2\pi),
$$
which is  a  circle in $S^2$.
For a positive integer $n\ge 2$, we set
$
\alpha=(\cos nt)dt.
$
Then, a flat front $f$ constructed from
the pair $(\alpha,\xi_{\phi})$ by Theorem~B is complete, since
$\cos nt$ is $L_2$-orthogonal to $\sin t$ and $\cos t$.
Since $\xi_{\phi}$ is a convex curve, $f$ 
has embedded ends.
All of the principal curvature lines are periodic
(cf. Proposition~\ref{prop:p-period}).
When $n=2$, there are four swallowtails and the
other singular points are cuspidal edges.
(See Figure~\ref{hamen}.)
The caustic $C_f$ is a cylinder over a closed
planar curves with $3/2$-cusps, 
which follows from the fact
that $\nu$ is a circle.
$C_f$ is co-orientable when $n$ is even and is non-co-orientable
when $n$ is odd.
\end{example}

We now prove Theorem~D
in the introduction. 

\medskip
\noindent
({\it Proof of Theorem~D.})
Let 
$$
\hat\xi=(\xi_1,\xi_2,\xi_3):S^1:=\R/2\pi\Z\to S^2
$$
be a convex spherical  curve,
that is, a simple closed regular curve in the unit sphere
with non-vanishing geodesic curvature. 
Let $f:S^1\times \R\to \R^3$ be a complete
flat front with embedded ends whose asymptotic
Gauss map is $\hat\xi(t)$.
Then there exists a function $a(t)\in C^\infty(S^1)$
such that
$$
\int_0^{2\pi}a(t)\hat\xi(t)dt=0 
$$
and the pair $(a(t)dt,\hat\xi)$ represents $f$ via the formula
\eqref{eq:B} in the introduction.

Suppose that $a(t)$ does not change sign.
We may assume that $a(t)>0$.
Since $\hat\xi(t)$ is a convex curve in $S^2$, we may also assume that
$\hat\xi(t)$ lies in the half-plane $D_+:=\{(x,y,z)\in \R^3\,;\, z>0\}$,
which yields that $a(t)\hat\xi(t)\in D_+$ and
the integral never vanishes, a contradiction.

Next, suppose that $a(t)$ changes sign exactly twice at 
$t=t_1,t_2$ ($0<t_1<t_2<2\pi$) where $S^1=\R/(2\pi \Z)$.
We define a subspace $\mathcal A$ of $C^\infty(S^1)$
by
$$
\mathcal A:=\{c_1 \xi_1+c_2 \xi_2+c_3 \xi_3
\in C^\infty(S^1)\;;\,
c_1,c_2,c_3\in \R \}.
$$
Since $\hat\xi(t)$ is a convex curve,
it meets each 2-dimensional linear subspace
$$
\left\{
(x,y,z)\in \R^3\,;\, c_1x+c_2y+c_3z=0,\,\, (c_1)^2+(c_2)^2+(c_3)^2
\ne 0\right\}
$$ at most twice.
Thus the linear map
$$
L:\mathcal A \ni \phi \mapsto (\phi(0),\phi(t_1),\phi(t_2))\in \R^3
$$
must be a linear isomorphism. 
So there exists a unique $\phi(t)\in \mathcal A$ such that
$$
(\phi(0),\phi(t_1),\phi(t_2))=(a(0),a(t_1),a(t_2))=(a(0),0,0).
$$
Since $\phi(t)$ changes sign at most two points, 
the sign of $\phi(t)$ coincides with 
that of $a(t)$, which implies that $a(t)\phi(t)>0$
holds for all $t\ne t_1,t_2$.
Thus we have  $\int_0^{2\pi}a(t)\hat\xi(t)dt>0$, a contradiction.
So $a(t)$ changes signs at least four times.

On the other hand, since $\hat\xi$ is a regular curve, $f$
has no linear singular points (see Proposition~\ref{prop:singular}).
Thus, a singular point which is not a cuspidal edge
contains in the zeros of $a(t)$ by (1) of 
Proposition~\ref{prop:singular}.
Thus there are four singular points which are not cuspidal edges.
\qed

\medskip
The conclusion of Theorem~D is optimal,
as the flat front with embedded end as in Figure~\ref{hamen} 
has exactly four swallowtails.
On the other hand, there exists a closed space curve
$c(t)$ with non-zero torsion function. (See Example~\ref{ex:twist} below.)
So we cannot drop the assumption of embeddedness of ends in Theorem~D. 
We also remark here that the above proof based on an analog
of the proof of the classical four vertex theorem for periodic functions
(see \cite{gmo} or \cite[Theorem A.4]{tu2} for details).
There are several refinements of the
classical four vertex theorem. See, for example,
\cite{ot}, \cite{u1} and  \cite{tu1}.
We give here the following application:

\begin{corollary}\label{cor:5-2}
Let $c:S^1\to \R^3$ be an embedded space curve
with non-vanishing curvature function $\kappa$.
Then the tangential developable 
$$
f(t,v):=c(t)+v c'(t)
$$
gives an admissible developable frontal.
Moreover, $f$ is a complete front if and only
if the torsion function $\tau$
never vanishes. 
In this case, $f$ has non-embedded ends.
\end{corollary}

\begin{proof}
Without loss of generality, we may assume that
$t$ is the arclength parameter of $c(t)$.
The unit bi-normal vector $\mathbf b(t)$ 
of the space curve $c(t)$ must be
the normal vector of $f$.
Since $f$ is complete, $f$ has no linear singularity.
Then by Theorem~\ref{thm:w-rep}, 
the fact that $f$ is a front implies that
$\tau$ has no zeros, since $\mb b'=\tau \mb n$, where
$\mb n:=c''(t)$.
In this case, as pointed out in Example~\ref{ex:dev},
$f$ is an admissible orientable developable front.
Moreover, by Proposition~\ref{prop:3-2}, it is
weakly complete. Then by Corollary~\ref{4-2},
$f$ is complete.
Suppose that $f$ has embedded ends.
Then Theorem~D yields that 
the curve $c(t)$ cannot be a regular curve,
a contradiction.
\end{proof}

\begin{figure}[htb]
  \begin{center}
         \includegraphics[width=4cm]{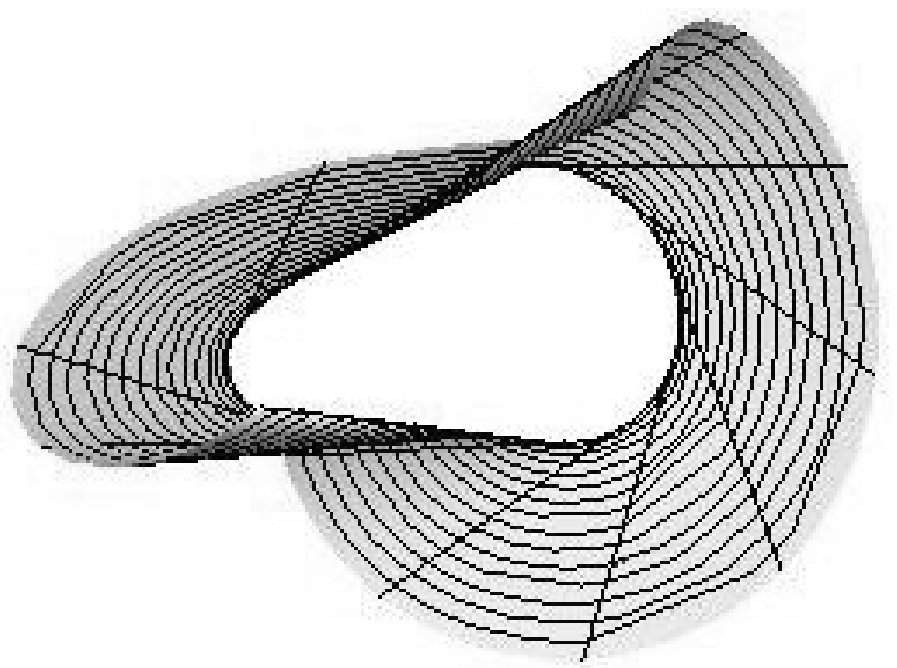}\qquad
         \includegraphics[width=4cm]{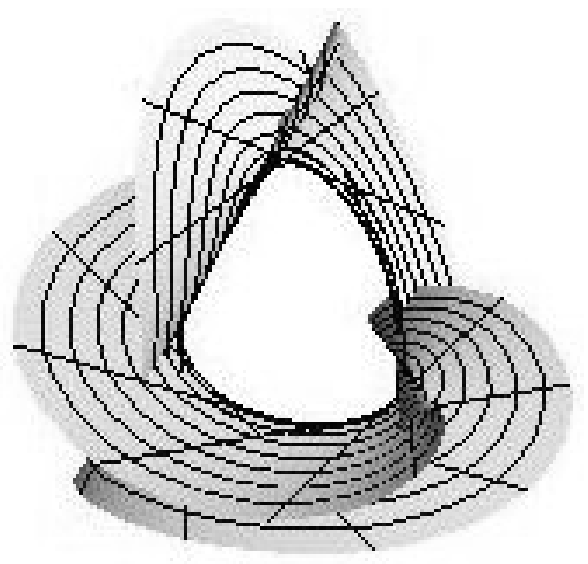}
\caption{$f(t,v)$ for $0\le v\le 7$ and 
$|v|\le 7$ in Example~\ref{ex:twist}}
\label{fig:tangential}
  \end{center}
\end{figure}

\begin{example}\label{ex:twist}
We set
$$
c(t):=((4+\cos 2t)\cos t,(4+\cos 2t)\sin t,\sin 2t)
\qquad (0\le t\le 2\pi),
$$
which gives a space curve of non-vanishing curvature
and of non-vanishing torsion.
By Corollary~\ref{cor:5-2}, the tangential developable
$f(t,v)$ of the curve gives a complete 
flat front with non-embedded ends (see Figure~\ref{fig:tangential} left).
Since $c$ has non-vanishing curvature function, 
all singular points are non-linear.
Moreover, since $c$ is a regular curve,
Proposition~\ref{prop:singular} yields that
the singular set $c([0,2\pi])$ consists of
cuspidal edges (see Figure~\ref{fig:tangential}, right).
\end{example}

%\newpage
\appendix

\section*{Appendix: Fronts in $S^2$ and their caustics.}

In this appendix, we shall define fronts in $S^2$
and explain their caustics.
(\cite{A} is a good reference.)
Let $J$ be an open interval.
A $C^\infty$-map $\gamma:J\to S^2$ is called
a {\it frontal}  if there exists
a  $C^\infty$-map $n:J\to S^2$ 
such that
$n(t)$ is perpendicular to $\gamma(t)$ and $\gamma'(t)$
for each $t\in J$.
The curve $\pm n(t)$ 
is called the {\it dual} of $\gamma(t)$.
By definition, $n(t)$ is also a frontal 
and $\gamma(t)$ is a dual of $n(t)$.
Moreover, since $\gamma'(t)$ and $n'(t)$ do not vanish
at the same time,
$\gamma(t)$ is called a {\it wave front} or a {\it front}.
Suppose now that $\gamma(t)$ is a front.
Then the  pair
$$
L=(\gamma(t),n(t)):J\to \{(x,y)\in S^2\times S^2 \,;\, x\cdot y=0\}=T^*S^2
$$
gives a Legendrian immersion with respect to 
the canonical contact structure, and the fronts can be
interpreted as projections of Legendrian immersions.
(Here $x\cdot y$ is the inner product of $x,y$ in $\R^3$.)
Now we set
$
\gamma_\theta(t):=\gamma(t)\cos \theta +  n(t)\sin \theta,
$
which is called the {\it parallel curve} of $\gamma$.
It can be directly checked that
$
n_\theta(t):=-\gamma(t)\sin \theta + n(t)\cos \theta  
$
gives the dual of $\gamma_\theta(t)$, that is,
{\it parallel curves of a front are also fronts}.
Since
$
n(t)=\gamma_{\pi/2}(t),
$
the dual front is in the parallel family of $\gamma$.

Now we assume that $\gamma(t)$ ($t\in J$)
is a regular curve.
Then we may assume that $t$ is the arc-length parameter 
of $\gamma$.
We set
$
n'(t)=-\kappa_g(t) \gamma'(t)
$,
then $\kappa_g(t)$ is the geodesic curvature of  $\gamma(t)$.
The parallel curve $\gamma_\theta(t)$ has
a singular point at $t=c$
if and only if
$\cos \theta - \kappa_g(c)\sin \theta=0$.
We define a smooth function
$
A(t):J\to \R/\pi \Z
$
such that
$$
\cos A(t) - \kappa_g\sin A(t)=0.
\leqno{(A.1)}
$$
Then
$
c(t):=\gamma(t)\cos A(t) +  n(t)\sin A(t)
$ 
is called the {\it caustic} (or {\it evolute}) of $\gamma(t)$, which 
has $\pm $ ambiguity, since $A(t)+\pi$ also satisfies (A.1).
The caustics $\pm c(t)$ give a 
parametrization of the singular set of the
parallel family of $\gamma$.
The following assertion 
was applied in Remark \ref{rem:duality}, and
can be proved directly.

\noindent
{\bf Proposition A.1.} \label{prop:A-2}
{\it Suppose that $\gamma(t)$ is
a regular spherical curve with arc-length parameter.
Then  $\gamma'$ is a regular spherical curve
and the dual of  $\gamma'$ coincides with
the caustic of $\gamma$.}

%%%%%%%%%%%%%%%%%%%%%%%%%%%%%%%%%%%%%%%%%%%%%%%%%%%%%%%%%%%%%%

\end{document}